\documentclass[12pt]{amsart}
\usepackage{amssymb}

\newtheorem{thm}{Theorem}[section]
\newtheorem{lem}[thm]{Lemma}
\newtheorem{cor}[thm]{Corollary}
\newtheorem{conj}[thm]{Conjecture}
\newtheorem{prop}[thm]{Proposition}
\newtheorem{question}[thm]{Question}

\theoremstyle{remark}

\newtheorem{rem}[thm]{Remark}

\theoremstyle{definition}
\newtheorem{defn}[thm]{Definition}

\numberwithin{equation}{section}

\allowdisplaybreaks[4]

\DeclareMathOperator{\rank}{rank}

\DeclareMathOperator{\Pic}{Pic}

\DeclareMathOperator{\Div}{Div}

\DeclareMathOperator{\supp}{supp}
\DeclareMathOperator{\Hom}{Hom}
\DeclareMathOperator{\Sym}{Sym}

\begin{document}

\vfuzz0.5pc
\hfuzz0.5pc 

\newcommand{\claimref}[1]{Claim \ref{#1}}
\newcommand{\thmref}[1]{Theorem \ref{#1}}
\newcommand{\propref}[1]{Proposition \ref{#1}}
\newcommand{\lemref}[1]{Lemma \ref{#1}}
\newcommand{\coref}[1]{Corollary \ref{#1}}
\newcommand{\remref}[1]{Remark \ref{#1}}
\newcommand{\conjref}[1]{Conjecture \ref{#1}}
\newcommand{\questionref}[1]{Question \ref{#1}}
\newcommand{\defnref}[1]{Definition \ref{#1}}
\newcommand{\secref}[1]{Sec. \ref{#1}}
\newcommand{\ssecref}[1]{\ref{#1}}
\newcommand{\sssecref}[1]{\ref{#1}}

\def \red{{\mathrm{red}}}
\def \tors{{\mathrm{tors}}}
\def \EQ{\Leftrightarrow}

\def \mapright#1{\smash{\mathop{\longrightarrow}\limits^{#1}}}
\def \mapleft#1{\smash{\mathop{\longleftarrow}\limits^{#1}}}
\def \mapdown#1{\Big\downarrow\rlap{$\vcenter{\hbox{$\scriptstyle#1$}}$}}
\def \smapdown#1{\downarrow\rlap{$\vcenter{\hbox{$\scriptstyle#1$}}$}}
\def \A{{\mathbb A}}
\def \I{{\mathcal I}}
\def \J{{\mathcal J}}
\def \CO{{\mathcal O}}
\def \C{{\mathcal C}}
\def \BC{{\mathbb C}}
\def \BQ{{\mathbb Q}}
\def \m{{\mathcal M}}
\def \H{{\mathcal H}}
\def \S{{\mathcal S}}
\def \Z{{\mathcal Z}}
\def \BZ{{\mathbb Z}}
\def \W{{\mathcal W}}
\def \Y{{\mathcal Y}}
\def \T{{\mathcal T}}
\def \P{{\mathbb P}}
\def \CP{{\mathcal P}}
\def \G{{\mathbb G}}
\def \F{{\mathbb F}}
\def \BR{{\mathbb R}}
\def \D{{\mathcal D}}
\def \L{{\mathcal L}}
\def \f{{\mathcal F}}

\def \closure#1{\overline{#1}}
\def \EQ{\Leftrightarrow}
\def \imply{\Rightarrow}
\def \isom{\cong}
\def \embed{\hookrightarrow}
\def \tensor{\mathop{\otimes}}
\def \wt#1{{\widetilde{#1}}}

\title{On Algebraic Hyperbolicity of Log Surfaces}

\author{Xi Chen}

\address{Department of Mathematics\\
South Hall, Room 6607\\
University of California\\
Santa Barbara, CA 93106}
\email{xichen@math.ucsb.edu}
\date{Mar. 30, 2001}
\thanks{Research supported in part by a Ky Fan Postdoctoral Fellowship}
\maketitle

\section{Introduction and Statement of Results}\label{s1}

\subsection{Introduction}

A complex manifold $M$ is hyperbolic in the sense of S. Kobayashi if
the hyperbolic pseudo-metric defined on $M$ is a metric. The study of
Kobayashi hyperbolicity (or just hyperbolicity in this paper) can be
roughly divided into two categories. One is the hyperbolicity of a
compact complex manifold. The other is the hyperbolicity of a compact
complex manifold with an ample divisor removed, or equivalently, the
hyperbolicity of an affine variety. In $\dim M = 1$, both problems are
completely settled. A compact Riemann surface $M$
(or a smooth projective curve) is hyperbolic if and only if its
geometric genus $g(M)$ is at least $2$. An affine algebraic curve $M =
C\backslash D$ is hyperbolic if and only if either $C\isom \P^1$ and
$\deg D \ge 3$ or $g(C) \ge 1$, where $C$ is a smooth projective
curve and $D$ is a reduced ample divisor on $C$. Equivalently, we may put
these two well-known facts into one sentence: $C\backslash D$ is
hyperbolic if and only if
\begin{equation}\label{s1:e0}
2 g(C) - 2 + \deg D > 0
\end{equation}
for a smooth projective curve $C$ and a reduced effective divisor
$D\subset C$.

Not much is known once the dimension goes beyond one,
even in the simplest case that $M$ is a hypersurface in $\P^n$ or $M$ is the
complement of a hypersurface in $\P^n$.
That is, we have the following famous Kobayashi conjecture. 

\begin{conj}[Kobayashi Conjecture]\label{conj1}
For $n\ge 2$,
\begin{enumerate}
\item a very general hypersurface $D\subset\P^{n+1}$ of
degree $\deg D\ge 2n + 1$ is hyperbolic;
\item $\P^n \backslash D$ is hyperbolic for a very general hypersurface
$D\subset \P^n$ of degree $\deg D \ge 2n + 1$.
\end{enumerate}
\end{conj}

We will refer the first part as Kobayashi's conjecture on hypersurfaces and
the second part as Kobayashi's conjecture on the complements of
hypersurfaces.

Most research on this conjecture was done for the case $n
= 2$, i.e., for surfaces in $\P^3$ and complements of plane curves in
$\P^2$. For surfaces in $\P^3$, a major step has been made
recently by J.P. Demailly and J. El Goul who proved that the conjecture is
true for $n = 2$ and $\deg D\ge 21$ \cite{DEG}. For complements of
plane curves in $\P^2$, Y.T. Siu and S.K. Yeung proved that the
conjecture is true for \(n = 2\) and \(\deg D \ge 10^{13}\) \cite{S-Y}.

Although hyperbolicity is essentially an analytic property,
it is closely related to the algebraic properties of the underlying
manifold, if the underlying manifold is an algebraic variety.
For example, there are no nonconstant maps from
an abelian variety to a hyperbolic algebraic variety.
Demailly proposed the following definition \cite[Chap. 2 \& 9]{D}
as an algebraic characterization of hyperbolicity.

\begin{defn}\label{defn1}
A smooth projective variety $X\subset\P^N$ is called {\it algebraically
hyperbolic\/} if there exists a positive real number $\epsilon$ such
that
\begin{equation}\label{defn1:e1}
2g(C) - 2\ge \epsilon \deg C
\end{equation}
for each reduced irreducible curve $C\subset X$, where $g(C)$ and
$\deg C$ are the geometric genus and the degree of the curve $C\subset\P^N$,
respectively. Obviously, the value of $\epsilon$ depends on the
polarization $X\embed \P^N$ of $X$, while being algebraically
hyperbolic is independent of the choice of polarization.
\end{defn}

Demailly proved that the following is true for a smooth projective
variety $X$
\begin{equation}\label{s1:e1}
\begin{split}
& X \text{ is hyperbolic} \imply X \text{ is algebraically hyperbolic}\\
&\imply \not\exists\text{ nonconstant maps from an abelian variety to }
X.
\end{split}
\end{equation}
Since algebraic hyperbolicity is a pure algebraic concept, it can be
handled entirely by algebro-geometrical methods. In general, it is
much easier to prove a variety is algebraically hyperbolic than to prove
it is hyperbolic. Thus \eqref{s1:e1} can be used to test which
varieties might be hyperbolic by checking whether they are
algebraically hyperbolic. On the other hand, from an
algebro-geometrical point of view, algebraic hyperbolicity is an
interesting notion itself, since finding the minimum genus of
a curve on an algebraic variety is an important question in algebraic
geometry \cite{Cl}.

It is a natural question to ask what is the similar notion of algebraic
hyperbolicity for an affine algebraic variety. We find it easier to
work with a compactification of an affine variety instead of
the variety itself. That is, we would like to work with a log pair
$(X, D)$ where $X$ is projective and $D$ is an effective
divisor. However, as we will see, the algebraic hyperbolicity we are
going to define for $X\backslash D$ is independent of the
compactification we choose.

It is obvious that if $X\backslash D$ is hyperbolic, there is no
nonconstant map from \(\BC^*\isom\P^1\backslash \{\text{2 points}\}\) to
$X\backslash D$, which is roughly equivalent to that each rational
curve $C\subset X$ meets $D$ at more than $2$ distinct points. This
suggests that the number of intersections between $C$ and $D$ should
appear on the LHS of \eqref{defn1:e1} for $(X, D)$. So it is natural
to define $(X, D)$ to be algebraically hyperbolic if there exists a
positive number $\epsilon$ such that
\begin{equation}\label{s1:e1.0}
2g(C) - 2 + |C\cap D| \ge \epsilon \deg C
\end{equation}
for all reduced irreducible curve $C\subset X$ with $C\not\subset
D$, where $|C\cap D|$ is the number of distinct points in
the set-theoretical intersection $C\cap D$. So if $(X, D)$ satisfies
\eqref{s1:e1.0}, every rational curve in $X$ should meet $D$ at no
less than three distinct points, which is what we want. However, there
is one problem with this definition. Ideally, we want our definition
of algebraic hyperbolicity of $(X, D)$ to be independent of the
compactifications of $X\backslash D$. However, $|C\cap D|$ does depend
on the choice of the compactifications of $X\backslash D$.
For example, suppose that $C$ meets $D$ at a simple triple point $p$
of $C$. Let $X'$ be the blowup of $X$ at $p$, $D'$ be the total
transform of $D$ and $C'$ be the proper transform of $C$. Then we have
$X\backslash D\isom X'\backslash D'$ but $|C'\cap D'| = |C\cap D| +
2$. Furthermore, if we assume that $C$ is rational and $p$ is the only
intersection between $C$ and $D$ in this example, then 
$(X, D)$ is not algebraically hyperbolicity by \eqref{s1:e1.0}.
However, the normalization of $C$ gives a morphism 
\(\P^1\backslash \{\text{3 points}\}\to X\backslash D\), which does
not violate the hyperbolicity of $X\backslash D$.
Therefore, we need first to have a better
notion of the number of intersections between $C$ and $D$.

\begin{defn}\label{defn2}
Let $P\subset D\subset X$, where $X$ is a quasi-projective
variety, $D\subset X$ is a closed subscheme of $X$ of codimension 1
and $P$ is a closed subscheme of $D$ of codimension 1. For each
reduced curve $C\subset X$ that meets $D$ properly, we define the
number $i_X(C, D, P)$ as follows. Let $\wt{X}$ be the blowup of $X$ at
$P$ and $\wt{C}$ and $\wt{D}$ be the corresponding proper transforms
of $C$ and $D$, respectively. Let $\nu: C^\nu \to \wt{C}\subset
\wt{X}$ be the normalization of $C$. Then $i_X(C, D, P)$ is the number
of distinct points in the set $\nu^{-1}(\wt{D})\subset C^\nu$. We
allow $P$ to be empty and write $i_X(C, D) = i_X(C, D, \emptyset)$.
Usually we abbreviate $i_X(C, D)$
and $i_X(C, D, P)$ to $i(C, D)$ and $i(C, D, P)$ if it is clear from
the context what the total space $X$ is.
Also let $i_{X, p}(C, D, P)$ be the local contribution of each point
$p$ to $i_X(C, D, P)$ such that $i_X(C, D, P) = \sum_{p\in X}
i_{X, p}(C, D, P)$ where we let $i_{X, p}(C, D, P) = 0$ if $p\not\in
C\cap D$.
\end{defn}

\begin{rem}\label{rem1}
Basically, we count the branches of $C$ at each point $p\in C\cap
D$ for $i_X(C, D)$. For example, $i_{X, p}(C, D) = 2$ if $C$ has a node
at $p\in C\cap D$ and $i_{X, p}(C, D) = 1$ if $C$ has a cusp at $p$.
So we always have $i_X(C, D)\ge |C\cap D|$. Sometimes we will call
$i_X(C, D)$ the number of intersections between $D$ and the
normalization of $C$.

It may look redundant to introduce the term $P$ in our definition
since we can simply use $i_{\wt{X}}(\wt{C}, \wt{D})$ for $i_X(C, D, P)$
after we blow up $X$ along $P$. Nevertheless, as we will see later,
$P$ arises in our problems in a natural way.
\end{rem}

Now we are ready to define algebraic hyperbolicity for log
varieties.

\begin{defn}\label{defn3}
Let $P\subset D\subset X$ be given as in \defnref{defn2} and
let $X\subset\P^N$ be projective. Then we call the triple $(X, D, P)$
{\it algebraically hyperbolic\/} if there exists a positive number
$\epsilon$ such that
\begin{equation}\label{defn3:e1}
2g(C) - 2 + i_X(C, D, P) \ge \epsilon \deg C
\end{equation}
for all reduced irreducible curves $C\subset X$ with $C\not\subset
D$. We call the log pair $(X, D)$ {\it algebraically hyperbolic\/} if
$(X, D, \emptyset)$ is algebraically hyperbolic.
\end{defn}

Here are a few (more or less) trivial comments.
\begin{enumerate}
\item Let $C$ be a smooth projective curve and $D\subset C$ be a reduced
effective divisor. Then $(C, D)$ is algebraically hyperbolic if and
only if \eqref{s1:e0} holds.
\item It is not hard to see that the algebraic hyperbolicity of $(X,
D)$ only depends on the
complement $X\backslash D$ and it is independent of the
compactification of $X\backslash D$. That is, if $(X, D)$ and $(X',
D')$ satisfy $X\backslash D\isom X'\backslash D'$, then $(X, D)$ is
algebraically hyperbolicity if and only if $(X', D')$ is. 
So instead of saying that $(X, D)$ is algebraically
hyperbolic, we may say $X\backslash D$ is algebraically hyperbolic.
\item If $(X_1, D_1)$ and $(X_2, D_2)$ are algebraically hyperbolic,
then the product
$(X_1, D_1)\times (X_2, D_2) = (X_1\times X_2, \pi_1^{-1}(D_1)
\cup \pi_2^{-1}(D_2))$ is also algebraically hyperbolic, where $\pi_1$ and
$\pi_2$ are the projections of $X_1\times X_2$ to $X_1$ and $X_2$,
respectively.
\item Being algebraically hyperbolic is a neither open nor close
condition. It is obviously not a close condition. To see it is not an
open condition either, take 
$(\P^1, \{\text{3 points}\})\times (\P^1, \{\text{3 points}\}) = (\P^1\times \P^1, D)$,
where $D$ is the union of three curves of type $(1, 0)$ and 
three curves of type
$(0,1)$. Such $(\P^1\times \P^1, D)$ is algebraically hyperbolic since
$(\P^1, \{\text{3 points}\})$ is. However, a generic deformation $D'$
of $D$ is an irreducible curve of type $(3, 3)$ and $(\P^1\times\P^1,
D')$ is not algebraically hyperbolic since there exists curves of type
$(1, 0)$ and $(0, 1)$ meeting $D'$ at only two distinct points.
\end{enumerate}

The term ``algebraically hyperbolic'' has been used before by several people
including the author in a weak sense. For example, in \cite{C3}, \cite{DSW}
and \cite{Ko}, ``algebraic hyperbolicity'' is used for the varieties $X$
that do not permit any nonconstant (algebraic) map from $C$ to $X$ with $C$ an
elliptic curve or 
\(C\isom \BC^{*}\isom \P^1\backslash \{\text{2 points}\}\). Obviously,
this notion of ``algebraically hyperbolic'' is weaker than what is
used here and hence we will call it ``weakly algebraically
hyperbolic'' to distinguish it from the definition used in this paper.

\begin{defn}\label{defn3.1}
A quasi-projective variety $X$ is called
{\it weakly algebraically hyperbolic\/} if it does
not contain any curve $C$ whose normalization is either an elliptic
curve or $\BC^{*}$.
\end{defn}

\subsection{Algebraic hyperbolicity of log surfaces}

Our first nontrivial example of an algebraically hyperbolic affine
variety is the complement of a very general plane curve $D\subset\P^2$
of degree at least five. Namely,
we will prove that Kobayashi's conjecture on the complements of
plane curves holds if we replace ``hyperbolic'' by ``algebraically
hyperbolic'', i.e., we have the following.

\begin{thm}\label{t1}
Let $D\subset \P^2$ be a very general plane curve of degree $d$.
Then
\begin{equation}\label{t1:e1}
2 g(C) - 2 + i_{\P^2}(C, D) \ge (d - 4) \deg C
\end{equation}
for all reduced irreducible curves $C\subset\P^2$ with $C\ne D$.
And hence $(\P^2, D)$ is algebraically hyperbolic if $\deg D = d\ge 5$.
\end{thm}

An immediate consequence of the above theorem is that \(\P^2\backslash
D\) is weakly algebraically hyperbolic for a very general quintic
curve \(D\), i.e., \(C\cap (\P^2\backslash D)\) is hyperbolic for all
curves \(C\subset\P^2\) with \(C\not\subset D\), or equivalently,
every holomorphic map from \(\BC\) to \(\P^2\backslash D\) is constant
if it is algebraically degenerated (see also \cite[Corollary, p. 612]{X2}).

It is a known fact that every rational curve on $\P^2$ meets a
very general quintic curve at no less than three distinct points (see
\cite{X2} and \cite{C3} for proofs and generalizations of this
statement). The above theorem may be regarded as another
generalization of this fact. It says that the number of the
intersections (counted in the way of \defnref{defn2}) between a
rational curve and a quintic curve grows at a rate proportional to the
degree of the rational curve. More precisely, it says the
following: fix a very general quintic curve $D\subset\P^2$ and then
for all nonconstant maps $f: C = \P^1\to\P^2$, there are at least
$\deg(f) + 2$ distinct points on $C$ mapping to points on $D$ by $f$.

Just as the original Kobayashi's conjecture, \thmref{t1} is much
easier to prove if $D$ is reducible. It was first proved by M. Green
that $\P^2\backslash D$ is hyperbolic for $D$ a general union of at
least four irreducible curves with total degree no less than $5$ \cite{G}.
Our next theorem has the same nature, i.e., it gives a criterion on
the algebraic hyperbolicity of $(X, D)$ when $D$ is reducible. Before
we state the theorem, let us first recall some definitions.

Let $X$ be a normal variety with canonical divisor $K_X$. We
call $X$ is {\it canonical\/} or has {\it canonical singularities\/}
if $K_X$ is $\BQ$-Cartier and for all proper birational morphisms $f:
Y\to X$, we have
\begin{equation}\label{s1:e1.1}
K_Y \sim_{num} f^* K_X + \sum a_i E_i
\end{equation}
with $a_i \ge 0$ for all exceptional divisors $E_i$, where
$\sim_{num}$ is the numerical equivalence. If \(X\) is a surface and
\(f: Y\to X\) is a minimal desingularization of \(X\), i.e.,
there are no \(f\)-contractible \(-1\) rational curves
\(E\subset Y\), then \(X\) is canonical iff \(K_Y\sim_{num}
f^* K_X\).

For a variety $X$, let $C_1(X)$ be the free abelian group generated by
the 1-cycles (curves) on $X$ and $N_1(X) = C_1(X)/\sim_{num}$.
We call $\varphi: N_1(X)\to \BR$ an additive function on $N_1(X)$ if
\(\varphi\in \Hom(N_1(X), \BR)\). Of course, if \(X\) is nonsingular,
\(\Hom(N_1(X), \BR) \isom N^1(X)\tensor \BR\) where \(N^1(X) =
\Div(X)/\sim_{num}\), i.e., every additive function \(\varphi\) on
\(N_1(X)\) is given by a \(\BR\)-divisor \(D\) such that \(\varphi(C)
= D\cdot C\) for all curves \(C\subset X\).

We call a closed subscheme $C\subset X$ {\it rigid\/} inside X
if there is no embedded deformation of $C$ inside $X$; otherwise, we
say $C$ is {\it nonrigid\/} inside $X$.

We call $\{I_1, I_2, ..., I_m\}$ a partition of a set $N$ if
$\cup_{k=1}^m I_k = N$ and $I_a\cap I_b = \emptyset$ for $1\le a\ne
b\le m$.

\begin{thm}\label{t2}
Let $S$ be a normal projective surface with canonical singularities,
$B = \sum_{i=0}^n B_i$ be an effective divisor of normal crossing on
$S$, $F\subset S$ be a curve on $S$ and \(P\subset B\cap F\).
Suppose that $B_i$ is a very general member of
a base point free (BPF) linear system \(\P \L_i\) for each $1\le i \le n$,
while \(B_0\) and \(F\) are fixed curves satisfying that
$B_0$ meets $F$ properly and $P\cap B_0$ is contained in the smooth locus
of $B_0$.

Let $\varphi: C_1(S) \to \BR$ be a function on
$C_1(S)$ and $\epsilon\in \BR$ be a real number satisfying
\begin{equation}\label{t2:e1}
(K_S + B - B_i)C \ge \epsilon \varphi(C)
\end{equation}
for $i =1,2,...,n$ and all nonrigid curves $C\subset S$,
\begin{equation}\label{t2:e2}
2g(C) - 2 + (B - B_0)C \ge \epsilon \varphi(C)
\end{equation}
for all reduced curves $C\subset S$ and
\begin{equation}\label{t2:e3}
2 g(C) - 2 + i_S(C, B, P) \ge \epsilon\varphi(C)
\end{equation}
for each reduced irreducible component $C\subset F$. Then
\begin{equation}\label{t2:e4}
2 g(C) - 2 + i_S(C, B, P) \ge \epsilon\varphi(C)
\end{equation}
for all reduced irreducible curves $C\subset S$ with
$C\not\subset B$.
\end{thm}

\begin{rem}\label{rem2}
Note that \(\varphi\) is an arbitrary function on \(C_1(S)\) (no
linearity is assumed). So what we are really saying in the above
theorem is that
\begin{equation}\label{rem2:e1}
\begin{split}
&\quad 2 g(C) - 2 + i_S(C, B, P)\\
&\ge
\min_{1\le i\le n}\left((K_S + B - B_i)C, 2g(C) - 2 + (B -
B_0)C\right)
\end{split}
\end{equation}
if \(C\not\subset F\) and \(C\) is not rigid. However, we choose to
present the theorem in the above form to make it consistent in
appearance with \thmref{t3}.
\end{rem}

It follows immediately from \thmref{t2} that
$(\P^2, \cup_{i=1}^d L_i)$ is algebraically hyperbolic for $d\ge 5$ lines
$L_1, L_2, ..., L_d$ in general position. To see this, let $n = d$,
$B_0 = \emptyset$, $B_i = L_i$, $\varphi(C) = \deg C$ and $P =
\emptyset$ and we have \eqref{t1:e1} with $D = \cup_{i=1}^d L_i$.
Now the question is how to go from this fact to \eqref{t1:e1} with $D$
irreducible. The natural approach is to degenerate an irreducible
plane curve to a union of lines. Note that this does not work for the
original Kobayashi's conjecture since being hyperbolic is
not an algebraic condition. However, this works for algebraic
hyperbolicity, although there are still some nontrivial technical
issues to be resolved. Actually, one can say that resolving these
technical issues is all this paper about. In summary, we can prove the
algebraic hyperbolicity of $(\P^2, D)$ by degenerating $D$ to a
union of $d$ lines. We will carry out this degeneration argument in a
more general setting as in the following theorem.

\begin{thm}\label{t3}
Let $S, B, B_i, \L_i, F$ be given as in \thmref{t2}.

Let $\{I_1, I_2, ..., I_m\}$ be a partition of $\{1, 2,..., n\}$,
$D_0 = B_0$ and $D_k$ be a very general member of the linear series
\begin{equation}\label{t3:e1}
\P \left(\bigotimes_{i\in I_k} \L_i\right)
\end{equation}
for $k = 1, 2, ..., m$. Let \(P\subset D\cap F\)
with \(P\cap D_0\) contained in the smooth locus of $D_0$, where 
\(D = D_0\cup D_1\cup D_2\cup ...\cup D_m\).

Let $\varphi: N_1(S)\to \BR$ be an additive function on
$N_1(S)$ and $\epsilon$ be a real number such that \eqref{t2:e1} and
\eqref{t2:e2} holds and
\begin{equation}\label{t3:e2}
2 g(C) - 2 + i_S(C, D, P) \ge \epsilon \varphi(C)
\end{equation}
for each reduced irreducible component $C\subset F$. Then
\begin{equation}\label{t3:e3}
2 g(C) - 2 + i_S(C, D, P) \ge \epsilon \varphi(C)
\end{equation}
for all reduced irreducible curves $C\subset S$
with $C\not\subset D$.
\end{thm}

It is obvious that what we are doing in the above theorem is
degenerating each curve $D_k$ to a union
\(\cup_{i\in I_k} B_i\) and hence $D$ to \(\cup B_i\), for which
\thmref{t2} can be applied.

To see how \thmref{t1} follows from \thmref{t3}, we take $n =
d$, $B_0 = \emptyset$, $B_1 = B_2 = ... = B_n$ to be the hyperplane
divisor of $\P^2$, $D \in \P H^0(\CO_S(B))$ to be a very general curve
of degree $d$, $\varphi(C) = \deg C$ and $P = \emptyset$;
obviously, we can choose $\epsilon = d - 4$ by \eqref{t2:e1} and
\eqref{t2:e2} and \thmref{t1} follows immediately. Also notice that
\(\{I_1, I_2, ..., I_m\}\) can be chosen to be an arbitrary partition of
\(\{1, 2, ..., n\}\) and consequently \eqref{t1:e1} holds for reducible
\(D\) as well. We have the following generalization of \thmref{t1} as a
corollary of \thmref{t3}.

\begin{cor}\label{cor0}
Let \(S\) be a normal projective surface with canonical singularities,
\(B\) be a BPF ample divisor on
\(S\) and \(D = \cup D_k\subset S\), where each \(D_k\) is a very general
member of \(\P H^0(\CO_S(d_k B))\) for some positive integer \(d_k\).
Let \(\alpha\) be the smallest number such that \(K_S + \alpha B\) is NEF.
Then
\begin{equation}\label{cor0:e1}
2 g(C) - 2 + i_S(C, D)\ge \left(d - \max(2, 1+\alpha)\right) B\cdot C
\end{equation}
for all reduced irreducible curves \(C\subset S\) with
\(C\not\subset D\), where \(d = \sum d_k\). Hence \((S, D)\) is
algebraically hyperbolic if \(d \ge 3\) and \(d\ge 2+\alpha\).

In particular,
\begin{equation}\label{cor0:e2}
2 g(C) - 2 + i_S(C, D)\ge \left(d - 4\right) B\cdot C
\end{equation}
for all reduced irreducible curves \(C\subset S\) with
\(C\not\subset D\) and
hence \((S, D)\) is algebraically hyperbolic if \(d \ge 5\).
In addition, if \((K_S + 2B)B\ge 0\), \eqref{cor0:e2} can
be improved by replacing \(d - 4\) on its RHS by \(d-3\).
\end{cor}

The second inequality \eqref{cor0:e2} is due to the fact that \(K_S +
3B\) is NEF, which is a consequence of Mori's cone theorem:
for a smooth projective variety \(X\) of dimension \(n\), the cone 
\(\overline{NE}(X)\) of effective curves of $X$ is generated by
\(K_X\)-NEF curves \(C\) (\(K_X\cdot C \ge 0\)) and the smooth
rational curves \(G\) with \(-(n+1)\le K_X\cdot G < 0\). More
generally, J. Koll\'ar proved that Mori's cone theorem holds for \(X\)
with isolated singularities which are locally the quotients of
isolated complete intersection singularities \cite{K}.
Obviously, surfaces with canonical singularities fall into this
category and hence Mori's cone theorem holds for canonical surfaces,
which is all we need in this paper. Actually, there is a much simpler
reason why Mori's cone theorem holds for a canonical surface
\(S\). Let \(f: \wt{S}\to S\) be a minimal desingularization of \(S\)
and then \(K_{\wt{S}} = f^* K_S\). Since Mori's cone theorem holds for
\(\wt{S}\), it is obvious that it holds for \(S\) as well.

If we take \(S = \P^2\) and \(D = \cup D_k\) in the
above corollary, where
\(D_k\subset \P^2\) is a very general plane curve of degree $d_k$, then
\begin{equation}\label{cor0:e3}
2 g(C) - 2 + i_{\P^2}(C, D) \ge (d - 4) \deg C
\end{equation}
for all reduced irreducible curves \(C\subset \P^2\) with
\(C\not\subset D\).
Hence \((\P^2, D)\) is algebraically hyperbolic if \(\deg D = d\ge 5\).
Note that it has been proved that such \(\P^2\backslash D\) is
hyperbolic if \(D\) has at least three components \cite{DSW}.

Similar statements to \thmref{t1} hold on other surfaces such as
rational ruled surfaces and Del Pezzo surfaces. Both can be derived from
\thmref{t3}.

\begin{cor}\label{cor1}
Let $\F_n$ be the rational ruled surface $\P (\CO\oplus \CO(n))$ over
$\P^1$ and let $M$ and $F$ be the divisors of $\F_n$ generating
$\Pic(\F_n)$ with $M^2 = -n$, $M\cdot F = 1$ and $F^2 = 0$.
Let \(D = \cup D_k\subset \F_n\), where each $D_k$ is a very
general member of a BPF complete linear series.
And let
\(D \in \P H^0(\CO_{\F_n}(aM + bF))\) with $a\le b$. Then
\begin{equation}\label{cor1:e1}
2 g(C) - 2 + i_{\F_n}(C, D) \ge \min\left(a - 3, b - an - 2\right) \deg C
\end{equation}
for all reduced irreducible curves $C\subset \F_n$ with
$C\not\subset D$, where $\deg C = (M + (n+1) F)C$.
And hence $(\F_n, D)$ is algebraically hyperbolic if
$b\ge a \ge 4$ and $b\ge 3+an$.
\end{cor}

\begin{cor}\label{cor2}
Let $\wt{\P}^2$ be the blowup of $\P^2$ at $2\le r\le 6$ general
points and let $R_1, R_2,..., R_n$ be all the $-1$
rational curves on $\wt{\P}^2$, i.e., $R_i$ is smooth and rational and
$K_{\wt{\P^2}}\cdot R_i = -1$ for $i=1,2,...,n$.
Let $D = \cup D_k\subset \wt{\P}^2$,
where each $D_k$ is a very general member
of a BPF complete linear series. Then
\begin{equation}\label{cor2:e1}
2 g(C) - 2 + i_{\wt{\P}^2}(C, D) \ge
\left(\min_{1\le i \le n} D\cdot R_i - 2\right) \deg C
\end{equation}
for all reduced irreducible curves $C\subset \wt{\P}^2$ with
$C\not\subset D$, where $\deg C = - K_{\wt{\P}^2}\cdot C$. And hence
$(\wt{\P}^2, D)$ is algebraically hyperbolic if $D \cdot R_i > 2$ for
all $R_i$.
\end{cor}

Naturally one may ask how sharp these results are. There is
no good answer to this question in general. However, in the cases we
have studied so far, our method seems to have produced sharp results. In
\thmref{t1}, it is not hard to see that the equality in \eqref{t1:e1}
is achieved when \(C\) is a bitangent or flex line of \(D\) so the
condition \(\deg D\ge 5\) is necessary for \((\P^2, D)\) to be
(algebraically) hyperbolic. In \coref{cor1}, if \(D\) is irreducible,
the equality \eqref{cor1:e1} is achieved when \(C = M\) or
\(C\) is a fiber of \(\F_n\to \P^1\) tangent to \(D\) so the conditions
\(b\ge a \ge 4\) and \(b\ge 3+an\) are necessary for \((\F^n, D)\) to
be (algebraically) hyperbolic. However, if \(D\) is reducible, we have
the exception that \(\P^1\times \P^1\backslash D\) is (algebraically)
hyperbolic for \(D\) a union of \(a = 3\) curves of type \((1, 0)\)
and \(b \ge 3\) curves of type \((0, 1)\). But with a little extra
effort, one can see that such \((\P^1\times \P^1, D)\)'s are the only
exceptions for which \coref{cor1} fails to be sharp. Finally, in
\coref{cor2}, the equality in \eqref{cor2:e1} is obviously achieved
by \(C = R_i\) that minimizes \(D\cdot R_i\) and the conditions
\(D\cdot R_i\ge 3\) are necessary for \((\wt{\P}^2, D)\) to be
(algebraically) hyperbolic.

So far we have been proving the algebraic hyperbolicity of the log
surface \((S, D)\) by specializing \(D\) while the
underlying surface $S$ is fixed ($S$ is rigid anyway in the cases we have
studied). However, there are situations that this approach does not
work. For example, let \((S, D)\)
be a log pair where $S$ is a very general K3 surface of genus $g$ and
$D$ is a very general curve in the primitive class of $S$. Since every
rational curve on $S$ meets $D$ at no less than $3$ distinct points if
$g\ge 3$, we expect $(S, D)$ to be algebraically hyperbolic 
if \(g\ge 3\). However, \thmref{t3} cannot be directly applied since
\(D\) is primitive. It turns out that we have to specialize the
underlying K3 surface $S$ at the same time we specialize $D$,
which results in the following theorem.

\begin{thm}\label{t4}
Let $S$ be a K3 surface of genus $g\ge 3$ and $D$ be a curve in the
primitive class of $S$. Then for a very general pair $(S, D)$,
\begin{equation}\label{t4:e1}
2 g(C) - 2 + i_S(C, D) \ge \frac{\lfloor (g - 3) /2\rfloor}{g - 1}
(C\cdot D)
\end{equation}
for all reduced irreducible curves $C\subset S$ with
$C\not\subset D$ and hence $(S, D)$ is algebraically hyperbolic if
$g\ge 5$.

Let $M$ be a general member of $\P H^0(\CO_S(n D))$ for some $n\ge 0$
and $P\subset D\cap M$. Then
\eqref{t4:e1} also holds if we replace $i_S(C, D)$ by $i_S(C, D, P)$.
\end{thm}

\begin{rem}\label{rem3}
Notice that
the above theorem fails to conclude the algebraic hyperbolicity of
$(S, D)$ for a K3 surface $S$ of genus $3$ and $4$, although we expect
such pairs to be algebraically hyperbolic.
The reason that we study such $(S, D)$ will be clear later.
\end{rem}

\subsection{Algebraic hyperbolicity of projective surfaces}

It has been speculated that there are connections between the
hyperbolicity of the hypersurfaces and that of their complements.
For example, it has been shown by Masuda and Noguchi
that for every positive integer $n$ there exists a number $d(n)$ such
that for each $d\ge d(n)$ there exists a smooth hypersurface 
\(D\subset \P^n\) of degree $d$ such that both $D$ and
$\P^n\backslash D$ are
hyperbolic \cite{Z}. However, the relation between the hyperbolicity
of $D$ and that of $\P^n\backslash D'$ is still quite mysterious, to
say the least. On the other hand, if we consider algebraic
hyperbolicity instead of hyperbolicity, we do have some partial answer
to the question. This forms the second part of this paper, which
studies the relation between the
algebraic hyperbolicity of projective surfaces, as defined by Demailly,
and the algebraic hyperbolicity of log surfaces, as defined
here. We will show that the study of algebraic hyperbolicity of
projective surfaces usually comes down to the study of that of log
surfaces via degeneration. For example, as we will see, the algebraic
hyperbolicity of a very general surface $S\subset \P^3$ of degree
$d\ge 6$ is more or less a consequence of the algebraic
hyperbolicity of $\P^2\backslash D$, where $D$ is a general union of
$d-1$ lines in $\P^2$. The link between these two objects is
provided by a notion of ``virtual genus'' of a curve lying on a
reducible variety.

\begin{defn}\label{defn4}
Let \(D = D_1\cup D_2\cup ...\cup D_n\) be a union of varieties
\(D_j\). For each $J\subset \{1,2,...,n\}$, let $D_J = \cap_{j\in J} D_j$ and
\(\partial D_J = D_J\cap (\cup_{k\not\in J} D_k)\).
Suppose that \(D\) is locally of normal crossing along
\(\partial D_j\) for each \(j\), i.e., for each \(J\subset
\{1,2,...,n\}\) with \(|J| \ge 2\), \(D\) is locally given by
\(\prod_{j\in J} x_j = 0\) at each point \(p\in D_J\) and \(p\not\in
\partial D_J\), where \(D_J = \emptyset\) if \(|J| > 1 + \dim D\).
Let $Q$ be a closed subscheme of $D$ of codimension $2$ satisfying
$Q\subset \cup (\partial D_j)$ and $\dim (Q\cap D_J) \le \dim D_J - 1$
for each $J\subset\{1,2,...,n\}$ with $|J| \ge 2$.
For each reduced irreducible curve $\Gamma\subset D$, let
$J\subset\{1,2,...,n\}$ be the subset such that $\Gamma\subset D_j$
for $j\in J$ and $\Gamma\not\subset D_j$ for $j\not\in J$ and we
define
\begin{equation}\label{defn4:e1}
\phi_Q(\Gamma) = 2 g(\Gamma) - 2 + i_{D_J}(\Gamma, \partial D_J, Q\cap
\partial D_J).
\end{equation}
Here if $\partial D_J = \emptyset$, we let
\(i_{D_J}(\Gamma, \partial D_J, Q\cap \partial D_J) = 0\).
Let $C\subset D$ be an arbitrary curve on \(D\) (possibly nonreduced
and reducible).
We use $\mu(\Gamma)$ to denote the multiplicity of an irreducible component
$\Gamma$ in $C$.
Then the {\it virtual genus\/} $g_Q^{vir}(C)$ of $C$ with respect to
$Q$ is defined by
\begin{equation}\label{defn4:e2}
2 g_Q^{vir}(C) - 2 = \sum_{\Gamma\subset C} \mu(\Gamma)\phi_Q(\Gamma)
\end{equation}
where we sum over all irreducible components \(\Gamma\subset C\).
\end{defn}

\begin{thm}\label{t5}
Let $X$ be a flat family of projective varieties over the disk $\Delta$
parameterized by $t$, where the
general fibers $X_t$ are irreducible and smooth and the central fiber
$X_0 = D = \cup_{j=1}^n D_j$ is locally of normal crossing 
along \(\partial D_j\) for each \(j\).
Let \(Q = \cup (\partial D_j) \cap X_{sing}\)
be the singular locus of
$X$ along \(\partial D_j\) and
suppose that \(\dim (Q\cap D_J)\le \dim D_J - 1\) for each
\(J\subset\{1,2,...,n\}\) with $|J| \ge 2$.
Let $Y$ be a flat family of curves over
$\Delta$ with the following commutative diagram:
\begin{equation}\label{t5:e1}
\begin{array}{ccc}
Y & \mapright{\pi} & X\\
\mapdown{} & & \mapdown{}\\
\Delta & \mapright{} & \Delta
\end{array}
\end{equation}
where $\pi: Y\to X$ is proper and $\Delta\to \Delta$ is a base change sending
$t$ to $t^\alpha$ for some $\alpha > 0$. Then 
\begin{equation}\label{t5:e2}
g(Y_t) \ge g_Q^{vir}(\pi(Y_0)),
\end{equation}
where we assume $Y$ to be reduced and
$g_Q^{vir}(\pi(Y_0))$ is the virtual genus of $\pi(Y_0)$ with
respect to $Q$, as defined in \defnref{defn4}.
\end{thm}

It is obvious by the above theorem that $X_t$ is algebraically
hyperbolic if $(D_J, \partial D_J, \partial D_J\cap Q)$ is algebraically
hyperbolic for each $J\subset \{1,2,...,n\}$. This is our basic
principle to reduce the algebraic hyperbolicity of
projective varieties to that of log varieties. Based on
this principle, we will study the algebraic hyperbolicity of
projective surfaces in \secref{s4}.
Our first theorem in this direction is the following, which concerns
the genus of a curve on a generic complete intersection in a projective
variety \(W\).

\begin{thm}\label{t6}
Let $W$ be a projective variety of dimension $n+2$ which is smooth in
codimension 2 and \(\P \L_1, \P \L_2, ..., \P \L_m\) be BPF linear systems
on $W$. Let \((a_{ij})_{m\times n}\) be a matrix of nonnegative
integers with \(\sum_{i=1}^m a_{ij} > 0\) for each $j$ and
\(\sum_{j=1}^n a_{ij} > 0\) for each $i$, 
let \(V_j\) be a very general member of the linear system
\begin{equation}\label{t6:e1}
\P \left(\bigotimes_{i=1}^m \L_i^{\tensor a_{ij}}\right)
\end{equation}
for \(j = 1,2,...,n\) and \(S = \cap_{j=1}^n V_j\).

Let \(\varphi: N_1(W)\to \BR\) be an additive function on \(N_1(W)\) and
\(\epsilon\) be a real number such that for all reduced curves
\(C\subset W\)
with \(C\cap W_{sing} = \emptyset\), the following holds:
\begin{equation}\label{t6:e2}
\left(K_W + \sum_{i, j} a_{ij} B_i - B_k\right) C \ge \epsilon \varphi(C)
\end{equation}
for $k = 1,2,...,m$ and
\begin{equation}\label{t6:e3}
2g(C) - 2 + \left(\sum_{i, j} a_{ij} B_i - \sum_{k=1}^m \mu_k B_k\right) C
\ge \epsilon \varphi(C)
\end{equation}
for all nonnegative integers \(\mu_k\) satisfying
\(\sum_{k=1}^m \mu_k = n\), \(\prod_{k=1}^m B_k^{\mu_k}\ne 0\) and
\(\mu_k \le \# \{a_{kj} \ne 0: 1\le j\le n\}\),
where $B_k$ is a general member of $\P\L_k$ for $k = 1,2,...,m$ and
\(\# \{a_{kj} \ne 0: 1\le j\le n\}\) is the number of nonzero entries
in the $k$-th row of the matrix \((a_{ij})\). Then
\begin{equation}\label{t6:e4}
2 g(C) - 2\ge \epsilon \varphi(C)
\end{equation}
for all reduced irreducible curves \(C\subset S\).
\end{thm}

Note that we assume very little on the total space \(W\) other than
the numerical properties, while most previous works on the problem
assume \(W\) to be homogeneous.

The most studied case of the above theorem is \(W = \P^3\),
\(m = n = 1\), \(\L_1 = H^0(\CO_{\P^3}(1))\) and
\(S\subset \P^3\) a very general surface of degree
\(d = a_{11}\). Here we take \(\varphi(C) = \deg C\). Then by the
above theorem, we may take \(\epsilon = d - 5\) and hence
\begin{equation}\label{t6:e5}
2 g(C) - 2 \ge (d - 5)\deg C
\end{equation}
for all reduced irreducible curves \(C\subset S\). This is a
well-known lower bound for $g(C)$ due to H. Clemens \cite{Cl}.
Clemens' result has been improved and generalized in various way by
L. Ein, C. Voisin, G. Xu, L. Chiantini, A. Lopez and Z. Ran  
(see \cite{E1}, \cite{E2}, \cite{V}, \cite{X1}, \cite{CLR} and
\cite{C-L}). G. Xu improved \eqref{t6:e5} by replacing
$\ge$ by $>$ in \cite{X1}, i.e.,
\begin{equation}\label{s1:e2}
g(C) > 1 + \frac{1}{2} (d - 5)\deg C.
\end{equation}
The significance of this improvement is that it settles Harris'
conjecture: a very general quintic surface \(S\subset \P^3\)
does not contain rational
or elliptic curves, or equivalently, $S$
is weakly algebraically hyperbolic in our language, while Clemens'
bound \eqref{t6:e5} only implies the nonexistence of rational curves
on $S$. Xu actually gave
a universal lower bound for \(g(C)\) independent of $\deg C$:
\begin{equation}\label{s1:e3}
g(C) \ge \frac{1}{2} d(d-3) - 2.
\end{equation}
Hence by \eqref{s1:e3}, a very general quintic surface does not have
any curve of genus $0$, $1$ and $2$. Xu's results \eqref{s1:e2} and
\eqref{s1:e3} was further generalized by Chiantini, Lopez and Ran in
\cite{CLR} and \cite{C-L}.

Xu's original proof of \eqref{s1:e2} is based on a rather subtle
deformation theory he developed for the pair $(S, C)$. L. Chiantini
and A. Lopez gave a shorter and more intrinsic proof of \eqref{s1:e2}
via the theory of {\it focal loci\/} in \cite{C-L} and their method
was generalized to deal with other surfaces. Later in the paper, we
will study two types of surfaces considered by them.

As another example of \thmref{t6},
take \(W = \P^{n_1}\times\P^{n_2}\times...\times\P^{n_m}\). Let \(B_1,
B_2, ..., B_m\) be the generators of \(\Pic W\), where \(B_i\) is the
pullback of the hyperplane divisor of \(\P^{n_i}\) under the
projection \(W\to \P^{n_i}\). And let \(\L_i = H^0(\CO_W(B_i))\) for
\(i=1,2,...,m\), \(V_j\in \P H^0(\CO_W(\sum_{i=1}^m a_{ij} B_i))\) and
\(S = \cap_{j=1}^n V_j\).
Set
\(\varphi(C) = \deg C = (\sum_{i=1}^m B_i) C\). Then by \thmref{t6},
we have
\begin{equation}\label{s1:e4}
2g(C) - 2 \ge \min_{1\le i\le m} \left(\sum_{j=1}^n a_{ij} - n_i - 2\right)
\deg C
\end{equation}
for all reduced irreducible curves \(C\subset S\). Therefore, \(S\)
is algebraically hyperbolic if \(\sum_{j=1}^n a_{ij} \ge n_i + 3\) for
\(i = 1,2,...,m\). It does not seem to have an easy way to derive
\eqref{s1:e4} using the methods of Clemens, Ein, Xu, etc.

Next, we will focus our attention to the complete intersections in
projective spaces. Let \(S\subset \P^{n+2}\) be a very general
complete intersection of type \((d_1, d_2, ..., d_n)\).
By \thmref{t6}, we have
\begin{equation}\label{s1:e5}
2 g(C) - 2 \ge \left(\sum_{j=1}^n d_j - n - 4\right) \deg C
\end{equation}
for all reduced irreducible curves $C\subset S$.
Hence $S$ is algebraically hyperbolic if 
\(d = \sum_{j=1}^n d_j > n+4\). This result \eqref{s1:e5} is
actually due to Ein. He generalized Clemens' theory on curves on
generic hypersurfaces to generic complete intersections in \cite{E1}
and \cite{E2}. Recently, Chiantini, Lopez and Ran further improved
Ein's result to the following:
\begin{equation}\label{s1:e6}
h^0\left(\CO_C(K_C - (d-n-4) f^*H )\right)\ge 1 + p
\end{equation}
for all maps \(f\) from a smooth curve \(C\) to \(S\),
where \(H\) is the hyperplane divisor of \(\P^{n+2}\) and \(p+1\) is
the dimension of the linear subspace of \(\P^{n+2}\) spanned by the
points on \(f(C)\). It is obvious that \(p\ge 1\) since \(S\) does not
contain any line. Therefore,
\begin{equation}\label{s1:e7}
2 g(C) - 2 \ge 2 + \left(d - n - 4\right) \deg C
\end{equation}
for all reduced irreducible curves \(C\subset S\). Take \(n = 1\) and
\eqref{s1:e7} becomes \eqref{s1:e2}.

It is possible to improve \eqref{s1:e5} along our line of argument. We
conjecture that the following is true:
\begin{equation}\label{s1:e8}
2 g(C) - 2 \ge \frac{2}{3}\prod_{j=1}^n d_j + 
\left(\sum_{j=1}^{n} d_j - n - 4\right) \deg C.
\end{equation}
However, the combinatorics involved is too daunting to be presented
here. Hopefully, we can treat this in a future paper.

Despite the improvement of \eqref{s1:e7} over \eqref{s1:e5}, it still
does not say much about the case that really interests us,
i.e., the algebraic hyperbolicity of complete intersections 
with \(\sum_{j=1}^{n} d_j = n+4\). Note that we do know such surfaces
are weakly algebraically hyperbolic by \eqref{s1:e7}.

\begin{question}\label{q1}
Is a very general complete intersection $S\subset \P^{n+2}$ of
type $(d_1, d_2, ..., d_{n})$ algebraically hyperbolic
if $\sum_{j=1}^{n} d_j = n+4$?  We may assume $S$ to be
nondegenerated, i.e., $d_j \ge 2$ for each $j$. Then $S$
is one of the following: a quintic surface in $\P^3$, a complete
intersection of type $(2, 4)$ or $(3, 3)$ in $\P^4$, a complete
intersection of type $(2, 2, 3)$ in $\P^5$ or a complete intersection
of type $(2, 2, 2, 2)$ in $\P^6$.
\end{question}

We will explain why we are interested in the above question
later.

Obviously, $n+4$ is the lower bound of $\sum d_j$ we can expect for a
complete intersection $S\subset\P^{n+2}$ of type 
\((d_1, d_2, ..., d_{n})\) to be hyperbolic since $S$ is rational
or K3 if \(\sum d_j < n+4\).

An analogous problem, though a little far-fetch, is Clemens'
conjecture on the finiteness of rational curves on a general quintic
threefold. Both of these two problems concern the cases where
the techniques of Clemens and Ein fail.

At present, we have a partial answer to the above question.

\begin{thm}\label{t7}
Let $S\subset \P^{n+2}$ be a very general complete intersection of
type \((d_1, d_2, ..., d_{n})\) with \(\sum_{j=1}^{n} d_j = n+4\).
Then $S$ is algebraically hyperbolic if $S$ is nondegenerated
and $n \ge 3$. More precisely, on such $S$,
\begin{equation}\label{t7:e1}
2 g(C) - 2 \ge \epsilon \deg C
\end{equation}
for each reduced irreducible curve $C\subset S$, where
\(\epsilon =1/6\) if $S\subset\P^5$ is of type $(2, 2, 3)$ and
\(\epsilon = 1/4\) if $S\subset\P^6$ is of type $(2, 2, 2, 2)$.
\end{thm}

We want to point out that our method of proving algebraic
hyperbolicity can be applied to a large class of projective surfaces,
not only generic complete intersections. Indeed, our approach is based upon
degenerations: as long as the projective surfaces under
consideration admit suitable degenerations, our method works.
For example, we have the following theorem concerning a certain class
of surfaces considered by Chiantini and Lopez in \cite{C-L}.

\begin{thm}\label{t8}
Let $D$ be a reduced irreducible curve in $\P^3$ and let $s, d$ be
two integers such that $d\ge s+4$ and
\begin{enumerate}
\item there exists a surface $R\subset\P^3$ of degree $s$ containing
$D$,
\item the general element of the linear system $\P
H^0(\CO_R(d)\tensor\CO_R(-D))$ is smooth and irreducible.
\end{enumerate}
Let $S\subset\P^3$ be a very general surface containing $D$. Then
\begin{equation}\label{t8:e1}
2g(C) - 2 \ge \epsilon \deg C
\end{equation}
for all reduced irreducible curve $C\subset S$ with
\begin{equation}\label{t8:e2}
\epsilon = \min\left(d - 5, d - s - 3,
\frac{2 g(D) - 2}{\deg D}, \frac{2 g(E) - 2}{\deg E}\right)
\end{equation}
where $E$ is a general member of 
\(\P H^0(\CO_R(d)\tensor\CO_R(-D))\), i.e., \(S\cap R = D\cup E\).
Therefore, $S$ is algebraically
hyperbolic if $d\ge 6$, $g(D)\ge 2$ and $g(E)\ge 2$.
\end{thm}

Of course, both \(g(E) = p_a(E)\) and \(\deg E\) are determined by
the given data \(d, s, D^2\) and \(\deg D\).

\begin{rem}\label{rem4}
The above theorem should be compared with Theorem 1.3 in \cite{C-L}, which
has the same hypotheses and concludes that
\begin{equation}\label{rem4:e1}
2 g(C) - 2\ge (d - s - 5)\deg C
\end{equation}
for all curves $C\subset S$ with $C\ne D, E$; therefore,
$S$ is algebraically hyperbolic if $d \ge s+6$, $g(D)\ge 2$ and
$g(E)\ge 2$.
There is a minor mistake in the original statement of this result in
\cite{C-L}, where it was stated
that \eqref{rem4:e1} holds for all $C\ne D$. The following simple example
shows that it is necessary to exclude $E$ from the consideration in
\eqref{rem4:e1}. Let \(s=1\), \(d \ge 6\) and
\(D\subset R\isom \P^2\) be a smooth and irreducible curve of degree
\(d - 1\). Obviously, \(E\) is a line in \(R\) and \eqref{rem4:e1}
fails if \(C = E\). So it is necessary to ask \(C\ne E\) in addition to
\(C\ne D\) in \eqref{rem4:e1}. However, the original proof in
\cite{C-L} still goes through after the change of the statement.
\end{rem}

In \cite{C-L}, Chiantini and Lopez also studied the algebraic
hyperbolicity of a very general projectively Cohen-Macaulay (PCM)
surface in \(\P^4\). A PCM surface \(S\subset\P^4\) is an irreducible
subvariety of \(\P^4\) of dimension two
whose ideal sheaf \(I_S\) has the minimal free resolution:
\begin{equation}\label{t9:e1}
0\rightarrow \bigoplus_{i=1}^{m} \CO_{\P^4}(-d_{2i})\rightarrow
\bigoplus_{j=1}^{m+1} \CO_{\P^4}(-d_{1j}) \rightarrow I_S\rightarrow 0.
\end{equation}
We call \((d_{11}, d_{12}, ..., d_{1,m+1}, d_{21}, d_{22}, ...,
d_{2m})\) the type of \(S\). Following the convention in
\cite{C-L}, we assume \(d_{11} \ge d_{12}\ge ...\ge d_{1,m+1}\) and 
\(d_{21}\ge d_{22}\ge ...\ge d_{2m}\) and let \(u_{ij} = d_{2i} -
d_{1j}\). It was shown in \cite{C-L} that
\begin{equation}\label{t9:e2}
2g(C) -2 \ge (u_{m,m+1} - 7)\deg C
\end{equation}
for all reduced irreducible curves \(C\subset S\). We will improve
\eqref{t9:e2} to the following.

\begin{thm}\label{t9}
Let \(S\) be a very general PCM surface in \(\P^4\) of type
\((d_{11}, d_{12}, ..., d_{1,m+1}, d_{21}, d_{22}, ...,
d_{2m})\). Suppose that \(u_{m,m+1} \ge 5\). Then
\begin{equation}\label{t9:e3}
2 g(C) - 2 \ge (u_{m,m+1} - 5) \deg C
\end{equation}
for all reduced irreducible curves \(C\subset S\) and hence \(S\) is
algebraically hyperbolic if \(u_{m, m+1}\ge 6\).
\end{thm}

\subsection{Some questions}

The first question we want to ask is whether a similar statement as
\eqref{s1:e1} holds for log varieties.

\begin{question}\label{q2}
Let $X$ be a smooth projective variety and let $D$ be an effective divisor
on $X$. Is it true that the hyperbolicity of $X\backslash D$ implies
the algebraic hyperbolicity of $(X, D)$?
\end{question}

We expect a positive answer to the above question although we are
unable to produce a proof at present. Let us analyze the proof of
\eqref{s1:e1} to see where the difficulty lies. This is taken
straight out of \cite[Theorem 2.1]{D}.

Let \(X\) be a polarized projective variety with the ample line bundle
\(L\), where \(X\) is
hyperbolic with the hyperbolic metric \(d s_X^2\). Let \(\omega_X\) be
the associated \((1,1)\) form of \(d s_X^2\) and \(\omega\) be a
positive \((1,1)\) form representing \(c_1(L)\). Since \(X\) is
compact, there exists a positive number \(\epsilon\) such that
\begin{equation}\label{q2:e1}
\omega_X \ge (\epsilon\pi) \omega.
\end{equation}
Let \(f: C\to X\) be a nonconstant morphism from a nonsingular curve \(C\)
of genus \(g(C) \ge 2\) to \(X\). Let \(d s_C^2\) be the hyperbolic
metric on \(C\) with constant Gaussian curvature \(-4\). Then by
Gauss-Bonnet,
\begin{equation}\label{q2:e2}
\int_C (-4) \frac{\omega_C}{2} = 2\pi c_1(T_C) = 2\pi \chi(C)
= 2\pi (2 - 2g(C))
\end{equation}
where \(\omega_C\) is the associated \((1,1)\) form of \(d
s_C^2\). A basic property of hyperbolic metrics is that
they are distant-decreasing under holomorphic maps, which means \(f^*
d s_X^2 \le d s_C^2\) in our case. This implies that
\begin{equation}\label{q2:e3}
\omega_C \ge f^*\omega_X.
\end{equation}
Combining \eqref{q2:e1}-\eqref{q2:e3}, we have
\begin{equation}\label{q2:e4}
\begin{split}
2 g(C) - 2 &= \frac{1}{\pi} \int_C \omega_C\\
&\ge \frac{1}{\pi} \int_C f^*\omega_X \ge 
\epsilon \int_C f^*\omega = \epsilon \deg f(C).
\end{split}
\end{equation}
Obviously, \eqref{q2:e1} is
the reason that the same argument cannot be directly carried over to a log
variety \((X, D)\): we may not have \eqref{q2:e1} for \(X\backslash
D\) since it is not compact.

Our second question is what remains of \questionref{q1}, especially
when \(S\) is a quintic surface.

\begin{question}\label{q3}
Is a very general quintic surface \(S\subset\P^3\) algebraically hyperbolic?
\end{question}

There is a famous conjecture in hyperbolic geometry by
M. Green and P. Griffiths, which says the following \cite{D}.

\begin{conj}[Green-Griffiths Conjecture]\label{conj2}
Let $X$ be a projective variety of general type.
Then every holomorphic map $f$ from $\BC$ to  $X$ is algebraically
degenerated, i.e., there exists a proper closed algebraic
subvariety \(Y\subset X\) such that
\begin{equation}\label{conj2:e1}
\bigcup_f f(\BC)\subset Y
\end{equation}
where \(f\) runs over all nonconstant holomorphic maps \(f:\BC\to X\).
\end{conj}

We call $f(\BC)$ an entire curve of $X$ for a nonconstant holomorphic map
\(f:\BC\to X\).

If this conjecture holds for surfaces,
every entire curve on a surface $S$ of general type is contained in a
rational or elliptic curve. The well-known theorem of Brody says that
a compact manifold is hyperbolic if and only if it does not contain
any entire curves \cite{B}. Thus if the conjecture of Green-Griffiths
holds for surfaces, a projective surface $S$ is hyperbolic
provided that $S$ does not contain any rational or elliptic curves,
i.e., $S$ is weakly algebraically hyperbolic. Since Xu has
proved that a very general quintic surface $S$ is weakly algebraically
hyperbolic, $S$ is hyperbolic if the conjecture of
Green-Griffiths holds. Actually, by \eqref{s1:e7},
all of the surfaces considered in \questionref{q1} are weakly
algebraically hyperbolic. Therefore, these surfaces are
hyperbolic if the conjecture of Green-Griffiths holds. So it
becomes an interesting question to ask whether such surfaces are
algebraically hyperbolic, since a negative answer to this question
will give counterexamples to the conjecture of Green-Griffiths.

The proof of \thmref{t7}, as we will see, comes down to the algebraic
hyperbolicity of a pair $(S, D)$, where $S$ is a K3 surface of genus
$5$ and $D$ is a curve in the primitive class of $S$. Using the
same argument for quintic surfaces, we see that the algebraic
hyperbolicity of a quintic surface follows 
from the algebraic hyperbolicity of $(S, D)$, where $S$ is a quartic K3
and $D$ is a curve in the primitive class of $S$.

Demailly asked the same question as \ref{q3} in \cite{D}.

Finally, we are looking for a generalization of the technique in this
paper to high dimensions. Specifically, we ask

\begin{question}\label{q4}
Let $D$ be a very general hypersurface of degree $d\ge 2n+1$ in
$\P^n$. Is $(\P^n, D)$ algebraically hyperbolic?
\end{question}

\subsection{Conventions}

\begin{enumerate}
\item Throughout the paper, we will work exclusively over $\BC$.
\item By a variety $X$ being very general, we mean that $X$ lies on
the corresponding parameter space (Hilbert scheme or
moduli space) with countably many proper closed subschemes removed. So
the notion of being very general relies on the fact that the base
field $\BC$ we work with is uncountable.
\end{enumerate}

\section{Proof of \thmref{t2} and \ref{t3}}\label{s2}

We start with a deformation lemma.

\subsection{A deformation lemma}\label{s2:ss1}


\begin{lem}\label{lem1}
Let $S$ be a normal projective surface with canonical singularities
and $B$ be an effective divisor of normal crossing on $S$.
Let $P\subset B$ be a finite set of
points contained in the nonsingular locus of $B$. For each integer
$\delta$, let $W_\delta$ be the scheme parameterizing the curves
$C\subset S$ with $2 g(C) - 1 + i_S(C, B, P) - (K_S + B)C = \delta$. Then
\begin{equation}\label{lem1:e0}
\dim W_\delta \le \max(0, \delta).
\end{equation}
In particular, \(\dim W_\delta \le 0\) if
\(\delta\le 0\), i.e., there are only
countably many reduced irreducible curves $C$ satisfying
\begin{equation}\label{lem1:e1}
2g(C) - 2 + i_S(C, B, P) < (K_S + B)C
\end{equation}
where $K_S$ is the canonical divisor of $S$.

On the other hand, if \(W_\delta\) is nonempty, then
\(\dim W_\delta \ge \delta - g\), where \(g\) is the genus of a
general member of \(W_\delta\).
\end{lem}

\begin{proof}
Let $C\subset S$ be a general member of $W_\delta$ and
let $\nu: C^\nu\to C\subset S$ be the normalization of $C$.

Notice that if we blow up $S$ along $P$ and let $\wt{S}$ be the resulting
surface and $\wt{C}$ and $\wt{B}$ be the proper transforms of $C$ and
$B$, respectively, we have $g(C) = g(\wt{C})$, $i_S(C, B, P) =
i_\wt{S}(\wt{C}, \wt{B})$ and $(K_S + B) C = (K_\wt{S} + \wt{B})
\wt{C}$. So we can simply start with $P = \emptyset$.

Let $\wt{S}\to S$ be a desingularization of $S$, $\wt{C}$ be the
proper transform of $C$ and $\wt{B}$ be the reduced total transform of
$B$. Our assumptions on $(S, B)$ guarantee that \((K_S + B) C =
(K_{\wt{S}} + \wt{B}) \wt{C}\). Also observe that \(i_S(C, B) =
i_{\wt{S}}(\wt{C}, \wt{B})\) and $\wt{B}$ has normal crossing.
So we may assume that $S$ is smooth.

Furthermore, if $C$ meets $B$ at a singular point $p$ of $B$, there
exists a series of blowups of $S$ over $p$ such that $\wt{C}$ does not meet
$\wt{B}$ at any singular point over $p$ on the resulting surface
$\wt{S}$. Therefore, we may assume that $C$ meets $B$ only at smooth
points of $B$. 

It suffices to show that the dimension of the deformation space of the
map $\nu: C^\nu\to S$ with $i_S(C, B)$ fixed is bounded from the above
by $\max(0, \delta)$. Here by the deformation of the map $\nu$, we
mean the deformation of the pair $(\nu, C^\nu)$.

Let $N$ be the normal bundle of the map $\nu: C^\nu\to S$. Let
$\nu^*(B) = m_1 p_1 + m_2 p_2 + ... + m_i p_i$ where $i = i_S(C, B)$.
Fixing $i_S(C, B)$ is equivalent to imposing certain
tangency conditions at $p_1, p_2, ..., p_i$.

The dimension of the deformation space of such map $\nu$ is at most
$h^0(C^\nu, N/N_\tors\tensor \CO(\sum_{k=1}^i (1 - m_k) p_k))$ (see
e.g. \cite[Chap. 2, Sec. B]{H}), where $N_\tors$ is the
torsion part of $N$. Since
\begin{equation}\label{lem1:e2}
\begin{split}
&\quad\deg \left(N/N_\tors\tensor
\CO\left(\sum_{k=1}^i (1 - m_k) p_k\right)\right)\\
&\le -(K_S + B) C + 2g(C) - 2 + i_S(C, B) = \delta - 1,
\end{split}
\end{equation}
\(h^0(C^\nu, N/N_\tors\tensor \CO(\sum_{k=1}^i (1 -
m_k) p_k)) \le \max(0, \delta)\) and we are done.

If \(W_\delta\) is nonempty, then \(W_\delta\) has at least the
virtual dimension
\begin{equation}\label{lem1:e3}
\begin{split}
& h^0\left(N\tensor \CO\left(\sum (1 - m_k) p_k\right)\right)\\
&- h^1\left(N\tensor \CO\left(\sum (1 - m_k) p_k\right)\right)
= \delta - g.
\end{split}
\end{equation}
\end{proof}

Now we are ready to prove \thmref{t2}.

\begin{proof}[Proof of \thmref{t2}]
There is nothing to do if $C\subset F$ since we have \eqref{t2:e3}.

Suppose that \(C\not\subset F\) is rigid in $S$. Since
$B_i$ moves in the BPF linear system \(\P \L_i\) for $i > 0$, by
B\'ezout, we can choose $B\backslash B_0$ such that
$C$ meets $B\backslash B_0$ transversely and
\((B\backslash B_0)\cap C\cap F = \emptyset\). Therefore,
\begin{equation}\label{t2:e5}
\begin{split}
2 g(C) - 2 + i_S(C, B\backslash B_0, P)&\ge 2 g(C) - 2 + (B - B_0)C\\
&\ge\epsilon\varphi(C).
\end{split}
\end{equation}

Suppose that $C\not\subset F$ is nonrigid. For each $1\le i \le n$, let us fix
$\cup_{j\ne i} B_j$. By \lemref{lem1}, there are only countably many
$C$ such that the following fails
\begin{equation}\label{t2:e6}
2 g(C) - 2 + i_S(C, \cup_{j\ne i} B_j, P) \ge (K_S + B - B_i) C \ge
\epsilon \varphi(C).
\end{equation}
Again, by B\'ezout, we can choose $B_i$ such that $B_i$ meets $C$
transversely, $B_i\cap B_j \cap C = \emptyset$ for each $j\ne i$ and
$B_i\cap C\cap F = \emptyset$ for all $C$ that \eqref{t2:e6} fails. 
In summary, we have either \eqref{t2:e6}
holds for some $i$, in which case we are done,
or $B\backslash B_0$ meets $C$ transversely and 
$(B\backslash B_0)\cap C\cap F = \emptyset$, in
which case we have \eqref{t2:e5}. Either way we obtain
\eqref{t2:e4}.
\end{proof}

\begin{rem}\label{rem6}
We take advantage of the fact that $B$ is reducible in the above
proof. By varying some component of $B$ while fixing the rest, we are
able to produce a family of curves for which \lemref{lem1} applies.
Another trick to use \lemref{lem1} is to ``vary'' $P$ if $C$ passes
through some points of $P$ and $P$ is chosen ``generically'' in some
sense. Suppose that $C$ passes through $\{ p_1, p_2, ..., p_l
\}\subset P$. Very often in our application,
we are allowed to fix part of $P$ and let the rest vary, e.g.,
fix $P_f = \{ p_1, p_2, ..., p_k\}$ and
let $P_v = \{p_{k+1}, p_{k+2}, ..., p_l\}$ vary for some $k < l$.
Then we will correspondingly obtain a one-parameter family of curves $C$.
Apply the lemma to $(C, B, P_f)$ and we have
\begin{equation}\label{rem6:e1}
2 g(C) - 2 + i_S(C, B, P_f) \ge (K_S + B)C.
\end{equation}
Notice that
\begin{equation}\label{rem6:e2}
i_S(C, B, P) = i_S(C, B, P_f) 
- \sum_{p\in P_v} \left(i_{S, p}(C, B) - i_{S, p}(C, B, P_v)\right)
\end{equation}
So we obtain
\begin{equation}\label{rem6:e3}
\begin{split}
& \quad 2 g(C) - 2 + i_S(C, B, P)\\
& \ge (K_S + B)C - \sum_{p\in P_v} \left(i_{S, p}(C, B) - i_{S, p}(C,
B, P_v)\right)\\
& \ge (K_S + B)C - \sum_{j=k+1}^l m_j
\end{split}
\end{equation}
where $m_j$ are the multiplicities of $C$ at $p_j$ for $j = 1,2,...,l$.
We will need this trick later in the proof of \thmref{t7}.
\end{rem}

The rest of the section is devoted to the proof of \thmref{t3}.

\subsection{Degeneration of \(D\)}\label{s2:ss2}

Here we will prove the theorem for $P = \emptyset$.
Since the case $P\ne \emptyset$ requires very little change in the
following proof, it will be left for the readers to verify the
theorem when $P\ne \emptyset$.

As mentioned before,
the proof of \thmref{t3} is carried out via a degeneration of $D$ to
a union \(\cup B_i\). However in practice, we will degenerate one
component at a time and argue by induction on the number of components
of \(D\). By that we mean, assuming \(1\in I_1\), at the first step,
we will fix \(D_2, D_3, ..., D_m\) and degenerate \(D_1\) to a union of two
curves \(B_1 \cup D_1'\) with \(B_1\in \P \L_1\) and \(D_1'\) a general
member of the linear series
\begin{equation}\label{t3:e4}
\P\left(\bigotimes_{i\in I_1\backslash\{1\}} \L_i\right).
\end{equation}
We can do this repeatedly. Each time we degenerate one component of
\(D\) and this will increase the number of components of \(D\) by
one until \(D\) is degenerated to a union \(\cup B_i\), for which
\thmref{t2} can be applied.

To set this up, let \(X = S\times\Delta\) and
\(W = \sum_{k=0}^m W^{(k)}\subset X\) be an effective divisor on $X$ with
components $W^{(k)}$ satisfying
\begin{enumerate}
\item \(W_t^{(k)} = D_k\) for \(k\ne 1\) and all \(t\in\Delta\);
\item \(W^{(1)}\) is a pencil of curves in the linear series
\begin{equation}\label{t3:e5}
\P\left(\bigotimes_{i\in I_1} \L_i\right)
\end{equation}
where the central fiber is \(W_0^{(1)} = B_1\cup D_1'\) while the
general fibers \(W_t^{(1)}\) (\(t\ne 0\)) are general members of the
linear series. 
\end{enumerate}
So the central fiber of $W$ is
$W_0 = D_0\cup (B_1\cup D_1')\cup D_2 \cup ...\cup D_m$.

Let $Y$ be a flat family of curves over $\Delta$ with the commutative
diagram \eqref{t5:e1}. We assume that $Y$ is irreducible,
$\pi: Y_t\to X$ maps $Y_t$ birationally onto its image for $t\ne 0$
and $\pi(Y)$ meets $W$ properly in $X$.
Our goal is to prove
\begin{equation}\label{s2:ss2:e1}
2 g(Y_t) - 2 + i_X(\pi(Y_t), W) \ge \epsilon \varphi\left(\pi(Y_t)\right)
\end{equation}
for \(t\ne 0\). The induction hypothesis is
\begin{equation}\label{t3:e6}
2 g(C) - 2 + i_S(C, W_0) \ge \epsilon \varphi(C)
\end{equation}
for each reduced irreducible curve \(C\subset S\) with \(C\not\subset
W_0\).

We may further apply semistable reduction to $Y\to X$ to make it into
a family of semistable maps with marked points $Y_t\cap \pi^{-1}(W)$
on $Y_t$ for $t\ne 0$. More specifically, we want the following to hold:
\begin{enumerate}
\item $Y$ is smooth and $Y_0$ is nodal;
\item $Y_t\cap \pi^{-1}(W)$ extends to disjoint sections of the fiberation
$Y\to\Delta$, i.e., the flat limit $\lim_{t\to 0} (Y_t\cap
\pi^{-1}(W))$ consists of $i_X(\pi(Y_t), W)$ distinct points
lying in the nonsingular locus of $Y_0$;
\item $Y\to X$ is minimal with respect to these properties.
\end{enumerate}

For each component $\Gamma\subset Y_0$, we define
$\sigma(\Gamma)\subset\Gamma$ to be the finite subset of $\Gamma$
consisting of the marked points
$\Gamma\cap \lim_{t\to 0} (Y_t\cap \pi^{-1}(W))$ and the intersections
between $\Gamma$ and the components of $Y_0$ other than \(\Gamma\).
We make the following observation:
\begin{equation}\label{s2:ss2:e2}
2 g(Y_t) - 2 + i_X(\pi(Y_t), W) = \sum_{\Gamma\subset Y_0}
\left(2p_a(\Gamma) - 2 + |\sigma(\Gamma)|\right)
\end{equation}
where $\Gamma$ runs over all irreducible components of $Y_0$ and
$p_a(\Gamma)$ is the arithmetic genus of $\Gamma$. Basically, each
marked point appears exactly once in $\coprod \sigma(\Gamma)$ and
counts one for $i_X(\pi(Y_t, W)$; each intersection between two
components of $Y_0$ appears twice in $\coprod \sigma(\Gamma)$ and
counts two for $2 g(Y_t) - 2 = 2 p_a(Y_0) - 2$. Thus \eqref{s2:ss2:e2}
follows.

With \eqref{s2:ss2:e2} in mind, we see that in order to show
\eqref{s2:ss2:e1}, it suffices to show
\begin{equation}\label{s2:ss2:e3}
2p_a(\Gamma) - 2 + |\sigma(\Gamma)|\ge \epsilon
\varphi\left(\pi(\Gamma)\right)
\end{equation}
for each irreducible component $\Gamma\subset Y_0$. If $\Gamma$ is contractible
under $\pi$, then
\begin{equation}\label{s2:ss2:e4}
2p_a(\Gamma) - 2 + |\sigma(\Gamma)| \ge 0
\end{equation}
since we assume $Y$ is minimal under the semistable reduction. And
hence \eqref{s2:ss2:e3} holds for all contractible components $\Gamma$.

We claim that \eqref{s2:ss2:e3} also holds for components
$\Gamma\subset Y_0$ with $\pi(\Gamma)\not\subset W_0$. Let
$\Gamma\subset Y_0$ be such a component, $\Sigma$ be the
reduced image of $\Gamma$ under $\pi$ and $\alpha$ be the degree of
the map $\pi_\Gamma: \Gamma\to\Sigma$, where \(\pi_\Gamma\) is the
restriction of \(\pi\) to \(\Gamma\).

Our basic observation is
\begin{equation}\label{s2:ss2:e5}
\pi_\Gamma^{-1}\left(\Sigma\cap W_0\right) \subset \sigma(\Gamma).
\end{equation}
Indeed, for each point $p\in \Sigma\cap W_0$ and $q\in\pi_\Gamma^{-1}(p)$,
either there is an irreducible component $\Gamma' \ne \Gamma$ passing
through $q$ or $q$ is one of the marked points; either way, we have
$q\in \sigma(\Gamma)$ and hence \eqref{s2:ss2:e5} follows.

Let $\nu_\gamma: \Gamma^\nu\to\Gamma$ and
$\nu_\sigma: \Sigma^\nu \to \Sigma$ be the normalizations of $\Gamma$
and $\Sigma$, respectively. We have the following commutative diagram
\begin{equation}\label{t3:e7}
\begin{array}{ccc}
\Gamma^\nu & \mapright{\pi_\Gamma} & \Sigma^\nu\\
\mapdown{\nu_\gamma} & & \mapdown{\nu_\sigma}\\
\Gamma & \mapright{\pi_\Gamma} & \Sigma
\end{array}
\end{equation}
where we use $\pi_\Gamma$ to denote both the maps
\(\Gamma^\nu\to\Sigma^\nu\) and
\(\Gamma\to\Sigma\); hopefully, this will not cause any confusion.

Then by \eqref{s2:ss2:e5} and observing that
$\sigma(\Gamma)$ is contained in the smooth locus of $\Gamma$, we have
\begin{equation}\label{t3:e8}
(\nu_\sigma\circ \pi_\Gamma)^{-1}(\Sigma\cap W_0)\subset
\nu_\gamma^{-1}(\sigma(\Gamma)).
\end{equation}
Note that 
\(\left|\nu_\sigma^{-1}(\Sigma\cap W_0)\right| =
i_S(\Sigma, W_0)\) by the definition of \(i_S(\Sigma, W_0)\) and
\(\left|\nu_\gamma^{-1}(\sigma(\Gamma))\right| =
|\sigma(\Gamma)|\). Hence the total ramification index of the map 
\(\pi_\Gamma:\Gamma^\nu\to\Sigma^\nu\) over
\(\nu_\sigma^{-1}(\Sigma\cap W_0)\) is at least
\begin{equation}\label{t3:e9}
\alpha\left(i_S(\Sigma, W_0)\right) - |\sigma(\Gamma)|.
\end{equation}
Apply Riemann-Hurwitz to \(\pi_\Gamma:\Gamma^\nu\to\Sigma^\nu\) and we have
\begin{equation}\label{s2:ss2:e6}
2 g(\Gamma) - 2 \ge \alpha(2g(\Sigma) - 2) + \alpha\left(
i_{S}(\Sigma, W_0)\right) - |\sigma(\Gamma)|.
\end{equation}
And by the induction hypothesis \eqref{t3:e6}
\begin{equation}\label{t3:e10}
2g(\Sigma) - 2 + i_{S}(\Sigma, W_0)\ge \epsilon \varphi(\Sigma).
\end{equation}
Combine \eqref{s2:ss2:e6} and \eqref{t3:e10} and we obtain
\eqref{s2:ss2:e3} for \(\Gamma\subset Y_0\) with
\(\pi(\Gamma)\not\subset W_0\). The way of applying Riemann-Hurwitz as
above is quite typical of our proof and it will appear
repeatedly in our argument.

It remains to show that \eqref{s2:ss2:e3} holds for all components
$\Gamma$ with \(\supp \pi(\Gamma)\subset W_0\). We will argue case by
case in the following order:
\begin{enumerate}
\item \(\supp \pi(\Gamma)\subset D_0\);
\item \(\supp \pi(\Gamma)\subset D_1'\);
\item \(\supp \pi(\Gamma)\subset \cup_{k\ge 2} D_k\);
\item \(\supp \pi(\Gamma)\subset B_1\).
\end{enumerate}

Let \(\pi_\Gamma\) be the restriction of \(\pi\) to \(\Gamma\)
and let \(\alpha\) be the degree of the map \(\pi_\Gamma\).

Suppose that \(\Gamma\) dominates an irreducible component \(\Sigma
\subset D_0\). For \(k\ne 0\), \(\Sigma\) meets
\(W_0^{(k)}\) properly and hence
\begin{equation}\label{t3:e11}
\Gamma \not\subset \pi^{-1}\left(W^{(k)}\right) \text{ and }
\pi_\Gamma^{-1}\left(\Sigma \cap W_0^{(k)}\right) \subset 
\pi^{-1}\left(W^{(k)}\right).
\end{equation}
Consequently,
\begin{equation}\label{t3:e12}
\pi_\Gamma^{-1}\left(\Sigma \cap W_0^{(k)}\right) \subset
\sigma(\Gamma)
\end{equation}
for \(k\ne 0\).
Using Riemann-Hurwitz in the same way as in \eqref{s2:ss2:e6}, we
obtain
\begin{equation}\label{t3:e13}
\begin{split}
2 g(\Gamma) - 2 + \left|\sigma(\Gamma)\right|
& \ge \alpha\left(2g(\Sigma) - 2 + i_S(\Sigma, W_0\backslash D_0)\right)\\
&\ge \alpha\left(2g(\Sigma) - 2 + (B -
B_0)\Sigma\right) \ge \epsilon \alpha \varphi(\Sigma).
\end{split}
\end{equation}
Therefore, \eqref{s2:ss2:e3} holds for all components \(\Gamma\)
with \(\supp \pi(\Gamma) \subset D_0\).

Suppose that \(\Gamma\) dominates an irreducible component
\(\Sigma\subset D_1'\). For \(k \ne 1\), \(\Sigma\) meets
\(W_0^{(k)}\) properly and hence
\begin{equation}\label{t3:e14}
\Gamma \not\subset \pi^{-1}\left(W^{(k)}\right) \text{ and }
\pi_\Gamma^{-1}\left(\Sigma \cap W_0^{(k)}\right) \subset 
\pi^{-1}\left(W^{(k)}\right).
\end{equation}
Consequently,
\begin{equation}\label{t3:e15}
\pi_\Gamma^{-1}\left(\Sigma \cap W_0^{(k)}\right) \subset
\sigma(\Gamma)
\end{equation}
for \(k\ne 1\).
Using Riemann-Hurwitz in the same way as in \eqref{s2:ss2:e6}, we
obtain
\begin{equation}\label{t3:e16}
\begin{split}
2 g(\Gamma) - 2 + \left|\sigma(\Gamma)\right|
& \ge \alpha\left(2g(\Sigma) - 2 + i_S(\Sigma, W_0\backslash 
(B_1\cup D_1'))\right)\\
&\ge \alpha\left((K_S + D_1')\Sigma + (B - B_1 - D_1')\Sigma\right)\\
&= \alpha(K_S + B - B_1)\Sigma \ge \epsilon \alpha \varphi(\Sigma),
\end{split}
\end{equation}
where it may need some explanation for
\begin{equation}\label{t3:e17}
2 g(\Sigma) - 2 = (K_S + \Sigma)\Sigma = (K_S + D_1')\Sigma.
\end{equation}
Since \(D_1'\) is a general member of the BPF linear series
\eqref{t3:e4}, by B\'ezout, \(\Sigma\) is smooth and \(g(\Sigma) =
p_a(\Sigma)\). If \(D_1'\) is connected, then \(\Sigma = D_1'\) and
\eqref{t3:e17} is obvious; otherwise, if \(D_1'\) has two or more
connected components, \(\Sigma\cdot D_1' = \Sigma^2 = 0\) and
\eqref{t3:e17} follows too.
 
In summary, \eqref{s2:ss2:e3} holds for all components \(\Gamma\)
with \(\supp \pi(\Gamma) \subset D_1'\).

Suppose that \(\Gamma\) dominates an irreducible component
\(\Sigma\subset \cup_{k\ge 2} D_k\). Without the loss of the
generality, let us assume that \(\Sigma\subset D_2\).

For \(k \ne 2\), \(\Sigma\) meets
\(W_0^{(k)}\) properly and hence
\begin{equation}\label{t3:e18}
\Gamma \not\subset \pi^{-1}\left(W^{(k)}\right) \text{ and }
\pi_\Gamma^{-1}\left(\Sigma \cap W_0^{(k)}\right) \subset 
\pi^{-1}\left(W^{(k)}\right).
\end{equation}
Consequently,
\begin{equation}\label{t3:e19}
\pi_\Gamma^{-1}\left(\Sigma \cap W_0^{(k)}\right) \subset
\sigma(\Gamma)
\end{equation}
for \(k\ne 2\).
Using Riemann-Hurwitz in the same way as in \eqref{s2:ss2:e6}, we
obtain
\begin{equation}\label{t3:e20}
\begin{split}
2 g(\Gamma) - 2 + \left|\sigma(\Gamma)\right|
& \ge \alpha\left(2g(\Sigma) - 2 + i_S(\Sigma, W_0\backslash D_2)\right)\\
&\ge \alpha\left((K_S + D_2)\Sigma + (B - D_2)\Sigma\right)\\
&= \alpha(K_S + B)\Sigma \ge \epsilon \alpha \varphi(\Sigma),
\end{split}
\end{equation}
Therefore, \eqref{s2:ss2:e3} holds for all components \(\Gamma\)
with \(\supp \pi(\Gamma) \subset \cup_{k\ge 2} D_k\).

It remains to verify \eqref{s2:ss2:e3} for \(\Gamma\) with
\(\supp \pi(\Gamma) \subset B_1\). This turns out to be the hardest
case.

\subsection{Construction of a fan}\label{s2:ss3}

Suppose that $\Gamma$ dominates an irreducible component \(G\subset
B_1\).

Let $\wt{X}$ be the blowup of $X$ along
$G\subset X_0$ and let $R$ be the exceptional divisor. Obviously,
$\wt{X}_0 = R\cup S$, $R\cap S = G$ and $R$ is the ruled surface over
$G$ given by $\P N_{G/X}$, where $N_{G/X}$ is the normal bundle of $G$
in $X$. Obviously, $N_{G/X}$ fits in the exact sequence:
\begin{equation}\label{s2:ss3:e1}
0\mapright{} N_{G/X_0} \mapright{} N_{G/X} \mapright{} 
\left.N_{X_0/X}\right|_G
\mapright{} 0.
\end{equation}
The above exact sequence actually splits basically due to the fact
that $X$ is a trivial family over $\Delta$ (see \cite{C4}). Hence 
\(N_{G/X} =
\left.N_{X_0/X}\right|_G\oplus N_{G/X_0} = \CO_G \oplus \CO_G(G)\).

Such construction has been extensively used by Z. Ran in his works on
Severi varieties of plane curves \cite{R} and by C. Ciliberto and
R. Miranda in their works on Nagata conjecture \cite{CM1} and
\cite{CM2}. In their cases, \(X = \P^2\times \Delta\) and \(G\subset
X_0\) is a line in \(\P^2\). So \(\wt{X}_0\) is the union of \(\P^2\)
and \(R\isom \F_1\) meeting transversely along the line \(G\).
Ran calls such \(\wt{X}_0\) a {\it fan\/}.

Let $\wt{W}\subset \wt{X}$ be the proper transform of $W$ under the
blowup $\wt{X}\to X$ and correspondingly, let $\wt{W}^{(k)}$ be the proper
transform of \(W^{(k)}\) for \(k = 0, 1,2,...,m\).
Then the central fiber $\wt{W}_0$ of $\wt{W}$
consists of $E\cup (W_0\backslash G)$, where $E\subset R$ is
a curve satisfying
\begin{equation}\label{s2:ss3:e2}
E\cap G = (W_0\backslash G)\cap G.
\end{equation}
Moreover, \(f_*(E) = G\) by the projection \(f: R\to
G\); combining this fact and \eqref{s2:ss3:e2} determines the
numerical class of \(E\). More precisely, $E$ can be described as follows.

For each \(0\le i\le n\), we use the notation $\L_i^G$ to denote the
restriction of \(\L_i\) to \(G\). Here we let \(\P \L_0\) be the
linear system consisting only of \(B_0\).

We let $\L_I = \tensor_{i\in I} \L_i$ for each index set
$I\subset\{0,1,2,...,n\}$ and let \(\L_I^G\) be the restriction of
\(\L_I\) to \(G\).

Every curve \(C\) on $R = \P N_{G/X}$ naturally corresponds to a global
section of $\Sym^\alpha N_{G/X}^\vee \tensor L$ for some integer $\alpha$
and line bundle $L$ on $G$, where \(\alpha\) is the number such that
\(f_*(C) = \alpha G\) and we use $N_{G/X}^\vee$ to denote the dual of
$N_{G/X}$, i.e., the conormal bundle.
Actually, there is a natural identification
\begin{equation}\label{t3:e21}
H^0(\CO_R(\alpha G)\tensor f^* L) \isom
H^0(\Sym^\alpha N_{G/X}^\vee \tensor L)
\end{equation}
where \(\CO_R(G)\) is the tautological line bundle of \(R\).

With this identification, $E$ corresponds to a global section of
\begin{equation}\label{s2:ss3:e3}
N_{G/X}^\vee\tensor \CO_G(B) = \CO_G(B - B_1)\oplus \CO_G(B).
\end{equation}
Indeed, \(E\) corresponds to a section in the subspace of
\(H^0(\CO_G(B - B_1))\oplus H^0(\CO_G(B))\) given by
\begin{equation}\label{s2:ss3:e4}
\left(\L_{I_1\backslash\{1\}}^G \oplus \L_{I_1}^G\right)
\tensor B_0^G \tensor \bigotimes_{k=2}^m D_k^G
\end{equation}
where \(B_0^G\in \L_0^G\) is the restriction of \(B_0\) to \(G\) and
\(D_k^G\in \L_{I_k}^G\) is the restriction of \(D_k\) to \(G\).
Geometrically, this implies that \(E = G'\cup \f\),
where \(G'\) is given by a section in
\begin{equation}\label{s2:ss3:e5}
\L_{I_1\backslash\{1\}}^G \oplus \L_{I_1}^G
\end{equation}
and $\f$ is a union of \((B - \sum_{i\in I_1} B_i)G\) fibers of
\(f: R\to G\). It is easy to see that \(G'\) and \(\f\) satisfy
\begin{gather}\label{t3:e22}
G'\cap G = D_1'\cap G,\\ \label{t3:e23}
\f\cap G = (\cup_{k\ne 1} D_k)\cap G,\\ \label{t3:e24}
G' \subset \wt{W}_0^{(1)} \text{ and } \f\subset \cup_{k\ne 1}
\wt{W}_0^{(k)}.
\end{gather}
Due to our generic choices of
\((B_1, D_1', W^{(1)})\), \(G'\) corresponds to a general
section in the space \eqref{s2:ss3:e5}.
In particular, \(G'\) is irreducible and \(f: G'\to G\) is an
isomorphism.

The map \(\pi: Y\to X\) induces a rational map \(Y\to
\wt{X}\). After resolving the indeterminacy of this map, we have
\begin{equation}\label{s2:ss3:e6}
\begin{array}{ccc}
\wt{Y} & \mapright{\wt{\pi}} & \wt{X}\\
\mapdown{} & & \mapdown{}\\
Y & \mapright{\pi} & X.
\end{array}
\end{equation}
We may further apply semistable reduction to $\wt{Y}\to \wt{X}$ and
hence let us assume that $\wt{\pi}:
\wt{Y}\to \wt{X}$ is a family of semistable
maps with marked points $\wt{Y}_t\cap \wt{\pi}^{-1}(\wt{W})$ on
$\wt{Y}_t$ for $t\ne 0$.


Since $\Gamma\subset Y_0$ is noncontractible under $\pi$, there exists
a unique component $\wt{\Gamma}\subset\wt{Y}_0$ that maps birationally
onto $\Gamma$. Let \(\sigma(\wt{\Gamma})\subset \wt{\Gamma}\) be the
set of points on \(\wt{\Gamma}\) over the points in
\(\sigma(\Gamma)\) under the map \(\wt{\Gamma}\to\Gamma\). Since
\(\sigma(\Gamma)\) is contained in the smooth locus of \(\Gamma\),
\(|\sigma(\wt{\Gamma})| = |\sigma(\Gamma)|\). And it is obvious that
\(p_a(\Gamma) \ge p_a(\wt{\Gamma})\) so in order to prove
\eqref{s2:ss2:e3}, it is enough to show the following
\begin{equation}\label{s2:ss3:e8}
2p_a(\wt{\Gamma}) - 2 + \left|\sigma(\wt{\Gamma})\right|\ge
\epsilon\alpha \varphi(G)
\end{equation}
for each component $\wt{\Gamma}$ that dominates $G$ via the map
$f\circ \wt{\pi}: \wt{\Gamma}\to G$ of degree $\alpha$.


For such \(\wt{\Gamma}\), one of the following must be true:
\begin{enumerate}
\item \(\supp \wt{\pi}(\wt{\Gamma}) = G\);
\item \(\supp \wt{\pi}(\wt{\Gamma}) = G'\);
\item \(\wt{\pi}(\wt{\Gamma})\) meets \(G\cup E\) properly.
\end{enumerate}

If $\wt{\Gamma}$ dominates $G$, we claim that
\begin{equation}\label{s2:ss3:e9}
\wt{\pi}_\Gamma^{-1}\left(G\cap E\right) =
\wt{\pi}_\Gamma^{-1} \left(G\cap \wt{W}_0\right)
\subset \sigma(\wt{\Gamma})
\end{equation}
where \(\wt{\pi}_\Gamma: \wt{\Gamma} \to G\) is the restriction of
\(\wt{\pi}\) to \(\wt{\Gamma}\).
The argument for \eqref{s2:ss3:e9} is more delicate than our previous argument
for similar statements due to the way \(\sigma(\wt{\Gamma})\)
is defined. It is easy to see that
\begin{equation}\label{s2:ss3:e9.1}
\wt{\pi}_\Gamma^{-1}(G\cap G') \subset \sigma(\wt{\Gamma})
\end{equation}
using our previous argument. But we need to be a little more careful
when we justify
\begin{equation}\label{s2:ss3:e9.2}
\wt{\pi}_\Gamma^{-1}(G\cap \f) \subset \sigma(\wt{\Gamma}).
\end{equation}
Let \(p\in \wt{\Gamma}\) be a point on \(\wt{\Gamma}\) such that
\(\wt{\pi}(p) \in G\cap \f\). Let \(Z\) be the union of components
such that \(\wt{\Gamma}\cup Z\) is the preimage of \(\Gamma\)
under the map \(\wt{Y}_0\to Y_0\). Obviously, \(Z\) is contractible
under \(\wt{Y}_0\to Y_0\). If there are no components of \(Z\) passing
through \(p\), then it is obvious that \(p\in \sigma(\wt{\Gamma})\);
otherwise, let
\(Z'\) be the maximal connected component of \(Z\) that passes through
\(p\). If \(Z'\) is contractible under \(\wt{\pi}\),
it is again obvious that \(p\in \sigma(\wt{\Gamma})\);
otherwise, \(\wt{\pi}(Z')\) is
supported along \(M\subset R\) which is a fiber of \(f: R\to
G\). Obviously, \(M\) meets \(G'\) at a point. So either there is a
point \(q\in Z'\) such that
\begin{equation}\label{s2:ss3:e9.3}
q\in\lim_{t\to 0} \left(\wt{Y}_t\cap
\wt{\pi}^{-1}(\wt{W}^{(1)})\right)
\end{equation}
or \(Z'\) intersects a component \(N\subset \wt{Y}_0\) that
dominates \(G'\). Either way, we see that \(Z'\) contracts to a point
in \(\sigma(\Gamma)\) under the map \(\wt{Y}_0\to Y_0\) and hence
\(p\in \sigma(\wt{\Gamma})\). Thus \eqref{s2:ss3:e9.2} follows.

Applying Riemann-Hurwitz to \(\wt{\pi}_\Gamma: \wt{\Gamma} \to G\)
in the same way as in \eqref{s2:ss2:e6}, we obtain \eqref{s2:ss3:e8} as follows
\begin{equation}\label{s2:ss3:e10}
\begin{split}
&\quad 2 g(\wt{\Gamma}) - 2 + \left|\sigma(\wt{\Gamma})\right|\\
&\ge\alpha(2 g(G) - 2 + G\cdot E) = \alpha\left((K_S + G) G + (B -
G)G\right)\\
&= \alpha (K_S + B) G \ge \epsilon \alpha \varphi(G).
\end{split}
\end{equation}

Suppose that \(\wt{\Gamma}\) dominates \(G'\). We claim that
\begin{equation}\label{s2:ss3:e11}
\wt{\pi}_\Gamma^{-1}\left(G' \cap G\right) \subset \sigma(\wt{\Gamma}).
\end{equation}
To see this, we need the following local lemma.

\begin{lem}\label{lem2}
Let \(X\subset \Delta_{x_1 x_2 ... x_r}^{r}\times \Delta\) be the
hypersurface given by $x_1x_2...x_n = t$ for some $n \le r$,
where \(\Delta_{x_1 x_2 ... x_r}^{r}\)
is the $r$-dimensional polydisk parameterized by
\((x_1, x_2, ..., x_r)\) and $\Delta$ is the disk parameterized by $t$. 
Let \(X_0 = D = \cup_{k=1}^n D_k\) with \(D_k = \{ x_k = t = 0 \}\).

Let
$Y$ be a flat family of curves over $\Delta$ with the commutative diagram
\eqref{t5:e1}. Suppose that there exists a curve \(\Gamma\subset
\pi(Y_0)\) with \(\Gamma \subset D_i\) for some \(i\). 
Then for each \(D_k\) there exists a curve \(\Gamma' \subset
\pi(Y_0)\) with \(\Gamma' \subset D_k\).
\end{lem}

Before we prove \lemref{lem2}, let us see how to justify
\eqref{s2:ss3:e11} by the lemma.

At each point \(p\in G\cap G'\), \(\wt{X}\) is locally given by
\begin{equation}\label{s2:ss3:e12}
\Delta_{x_1 x_2 x_3}^3\times \Delta / (x_1x_2 = t).
\end{equation}
Let \(q \in \wt{\pi}_\Gamma^{-1}(p)\).
Then by the above lemma, there exists a component \(M\subset\wt{Y}_0\)
and a union \(Z\subset\wt{Y}_0\) of contractible components such that
\begin{enumerate}
\item \(M\) maps to a curve on \(S\) passing through \(p\);
\item \(Z\) is connected and \(Z\) contracts to the point \(p\) by
\(\wt{\pi}\);
\item \(Z\) meets \(\wt{\Gamma}\) at \(q\) and \(M\) at a point \(q'\)
that \(\wt{\pi}(q') = p\).
\end{enumerate}
Roughly, \(\wt{\Gamma}\) and \(M\) are joined by a connected
union of contractible components over \(p\). The same must be true for
\(\Gamma\) and \(M'\) on \(Y_0\), where \(M'\) is the image of \(M\)
under the map \(\wt{Y}_0\to Y_0\).
Therefore, \(q\in \sigma(\wt{\Gamma})\) and \eqref{s2:ss3:e11} follows.

\begin{proof}[Proof of \lemref{lem2}]
Without the loss of generality, let us assume that $Y$ is smooth and
irreducible. Suppose that it fails for $D_k$ for some \(k\). 
Since \(X\) is smooth,
\(D_k\) is a Cartier divisor of $X$ and hence \(\pi^*(D_k) = M\)
is a divisor of $Y$ which is supported on a union of contractible components of
$Y_0$. Then \(M^2 = M \cdot \pi^*(D_k) = 0\), which can happen only
when \(\supp M = Y_0\). 
But this contradicts with \(\Gamma\not\subset D_k\).
\end{proof}

\begin{rem}\label{rem5}
Some variations of the above lemma have been extensively used in \cite{C1},
\cite{C2} and \cite{C4}. It is also crucial in our proof of
\thmref{t5} in \ssecref{s4:ss3}.
\end{rem}

In addition to \eqref{s2:ss3:e11}, we also have
\begin{equation}\label{s2:ss3:e13}
\wt{\pi}_\Gamma^{-1}\left(G' \cap \f\right) = 
\wt{\pi}_\Gamma^{-1}\left(G' \cap (\cup_{k\ne 1}
\wt{W}_0^{(k)})\right)\subset \sigma(\wt{\Gamma}).
\end{equation}
Combining \eqref{s2:ss3:e11} and \eqref{s2:ss3:e13} and applying
Riemann-Hurwitz as usual, we derive \eqref{s2:ss3:e8} as follows
\begin{equation}\label{s2:ss3:e14}
\begin{split}
&\quad 2 g(\wt{\Gamma}) - 2 + \left|\sigma(\wt{\Gamma})\right|\\
&\ge\alpha(2 g(G') - 2 + G'\cdot G + G'\cdot \f)\\
&= \alpha\left((K_S + B_1) G + D_1'\cdot G + (B - B_1 - D_1') G
\right)\\
&= \alpha (K_S + B)G \ge \epsilon \alpha \varphi(G).
\end{split}
\end{equation}

It remains to prove \eqref{s2:ss3:e8} for \(\wt{\Gamma}\) with
\(\supp \wt{\pi}(\wt{\Gamma})\not\subset E\cup G\).
Let \(\wt{\Sigma}\) be the reduced image of the map 
\(\wt{\pi}: \wt{\Gamma}\to \wt{X}_0\) and \(\gamma\) be the degree of
the map \(\wt{\pi}: \wt{\Gamma}\to \wt{\Sigma}\). Then by our
``standard'' Riemann-Hurwitz argument again,
\begin{equation}\label{s2:ss3:e15}
2g(\wt{\Gamma}) - 2 + \left|\sigma(\wt{\Gamma})\right| \ge
\gamma\left(2g(\wt{\Sigma}) - 2 + i_R(\wt{\Sigma}, G\cup E)\right).
\end{equation}
Therefore, \eqref{s2:ss3:e8} follows from
\eqref{s2:ss3:e15} provided that we can prove
\begin{equation}\label{s2:ss3:e16}
2g(\wt{\Sigma}) - 2 + i_R(\wt{\Sigma}, G\cup E)\ge \beta
(\epsilon\varphi(G))
\end{equation}
where \(\beta = \alpha/\gamma\) is the degree of the map 
\(f: \wt{\Sigma}\to G\).

In summary, we finally reduce \thmref{t3} to \eqref{s2:ss3:e16}, which
is a statement concerning the log surface \((R, G\cup E)\). We can
translate \eqref{s2:ss3:e16} into the following proposition.

\begin{prop}\label{prop1}
Let \(G\) be a smooth projective curve, \(\P \L_1, \P \L_2, ..., \P \L_l\)
be BPF linear systems on \(G\) with general members \(B_k\in \P \L_k\)
(\(1\le k\le l\)), and \(R\) be the ruled surface over \(G\) given by
\(\P(\CO_G\oplus \CO_G(B_1))\). We will abuse the notation a little
bit by using $G$ to denote the section of $f: R\to G$ with 
\(G^2 = -\deg B_1\), i.e., \(\CO_R(G)\) gives the tautological bundle of
\(\P(\CO_G\oplus \CO_G(B_1))\). 

Let \(G'\subset R\) be a general member of the linear series
\begin{equation}\label{prop1:e1}
\begin{split}
\P \left( \tensor_{i=2}^l \L_i \oplus \tensor_{k=1}^l \L_k\right)
&\subset \P H^0(\CO_G(B - B_1) \oplus \CO_G(B))\\
&\isom \P H^0(\CO_R(G + f^* B))
\end{split}
\end{equation}
where \(B = \sum B_k\).
And let \(\f = \f_1\cup \f_2\cup ...\cup \f_c\) be a union of \(c\) distinct 
fibers of the projection \(f: R\to G\).

Then for each reduced irreducible curve $\Sigma\subset R$ that
dominates $G$ under the projection $f$ and meets $G'\cup \f\cup G$
properly, we have
\begin{equation}\label{prop1:e2}
\begin{split}
&\quad 2g(\Sigma) - 2 + i_R(\Sigma, G'\cup \f\cup G)\\
&\ge \beta \left(2g(G) - 2 + \deg B - \deg B_1 - \max_{2\le k\le l}\deg B_k
+ c\right)
\end{split}
\end{equation}
where \(\beta\) is the degree of the map \(f: \Sigma\to G\).
\end{prop}

Applying the above proposition to \((R, G\cup E) = (R, G'\cup \f\cup
G)\) yields
\begin{equation}\label{s2:ss3:e17}
\begin{split}
&\quad 2g(\wt{\Sigma}) - 2 + i_R(\wt{\Sigma}, G\cup E)\\
&\ge \beta\left(2g(G) - 2 + (B - B_1) G - \max_{i\in I_1\backslash
\{1\}} B_i \cdot G\right)\\
&=\beta\left((K_S + B_1) G + (B - B_1) G - \max_{i\in I_1\backslash
\{1\}} B_i \cdot G\right)\\
&= \beta \min_{i\in I_1\backslash\{1\}}\left(K_S + B -  B_i \right) G.
\end{split}
\end{equation}
And hence \eqref{s2:ss3:e16} follows. So it remains to prove
\propref{prop1}.

\subsection{Proof of \propref{prop1}}\label{s2:ss4}

\begin{proof}[Proof of \propref{prop1}]
We argue by induction on $l$.

When \(l = 1\), we apply Riemann-Hurwitz directly to the projection 
\(f: \Sigma\to G\) by observing that 
\begin{equation}\label{prop1:e3}
f^{-1}(G\cap \f) \subset \Sigma\cap \f.
\end{equation}
This yields
\begin{equation}\label{s2:ss4:e1}
2g(\Sigma) - 2 + i_R(\Sigma, \f) \ge \beta(2 g(G) - 2 + c)
\end{equation}
and hence \eqref{prop1:e2} follows.

For $l > 1$, we will use degeneration to bring down $l$. Since
$G'$ is chosen to be generic in the linear series given in
\eqref{prop1:e1}, it can be degenerated to a union 
\(G''\cup \f'\). Here \(\f' = \f_{c+1}\cup \f_{c+2}\cup ...\cup
\f_{c+b_l}\) is the union of \(b_l = \deg B_l\) fibers of \(R\to G\) over the
vanishing locus of \(B_l\in \P \L_l\) and \(G''\) is
a general member of the linear series
\begin{equation}\label{prop1:e4}
\begin{split}
\P \left( \tensor_{i=2}^{l-1} \L_i \oplus \tensor_{k=1}^{l-1} \L_k\right)
&\subset \P H^0(\CO_G(B' - B_1) \oplus \CO_G(B'))\\
&\isom \P H^0(\CO_R(G + f^* B'))
\end{split}
\end{equation}
where \(B' = B - B_l\). Our induction hypothesis is
\begin{equation}\label{prop1:e5}
\begin{split}
&\quad 2g(\Sigma) - 2 + i_R(\Sigma, G''\cup (\f'\cup \f)\cup G)\\
&\ge \beta \bigg(2g(G) - 2 + \deg B' - \deg B_1\\
&\quad \quad \quad \quad - \max_{2\le k\le l-1}\deg B_k
+ c + \deg B_l\bigg)\\
&= \beta \left(2g(G) - 2 + \deg B - \deg B_1 - \max_{2\le k\le l-1}\deg B_k
+ c\right)
\end{split}
\end{equation}
if \(\Sigma\) meets \(G''\cup (\f'\cup \f)\cup G\) properly.

Let $X = R\times\Delta$ and $W = W^{(1)} + W^{(2)} + W^{(3)}$ be an effective
divisor on $X$ with components $W^{(k)}$ ($1\le k\le 3$) satisfying
\begin{enumerate}
\item \(W_0^{(1)} = G''\cup \f'\);
\item \(W_t^{(1)}\) is a general member of the linear series
\eqref{prop1:e1} for \(t\ne 0\);
\item \(W_t^{(2)} = \f\) and \(W_t^{(3)} = G\) for all \(t\in
\Delta\).
\end{enumerate}
So \(W_0 = G'' \cup \f' \cup \f \cup G\).

Let $Y$ be a flat family of curves over $\Delta$ with the commutative
diagram \eqref{t5:e1}. We assume that $Y$ is smooth and irreducible,
the map \(\pi: Y_t\to X\) is birational onto its image for $t\ne
0$ and \(\pi(Y)\) meets $W$ properly in $X$. Our goal is to prove
\begin{equation}\label{prop1:e6}
\begin{split}
&\quad 2 g(Y_t) - 2 + i_X(\pi(Y_t), W)\\
&\ge \beta \left(2g(G) - 2 + \deg B - \deg B_1 - \max_{2\le k\le l}\deg B_k
+ c\right)
\end{split}
\end{equation}
for $t\ne 0$, where $\beta$ is the degree of the map $f\circ\pi: Y_t\to G$.

As before, we may assume that $Y\to X$ is a family of semistable maps with
marked points $Y_t\cap \pi^{-1}(W)$ on $Y_t$ ($t\ne 0$) and
let \(\sigma(\Gamma)\subset\Gamma\) 
be defined similarly for each component $\Gamma\subset Y_0$.
Since we still have \eqref{s2:ss2:e2}, the proof of
\eqref{prop1:e6}, like that of \eqref{s2:ss2:e1},
again comes down to the estimation of \(2 p_a(\Gamma) - 2 +
\left|\sigma(\Gamma)\right|\) for each irreducible component
\(\Gamma\subset Y_0\).

Let \(\gamma\) be the degree of the map \(f\circ \pi: \Gamma\to G\).

If \(\supp \pi(\Gamma) = G\), then
\begin{equation}\label{prop1:e7}
\pi_\Gamma^{-1}\left(G\cap \left(G'' \cup \f'\cup \f\right)\right)
= \pi_\Gamma^{-1}\left(G\cap \left(W_0^{(1)} \cup W_0^{(2)}\right)\right)
\subset \sigma(\Gamma)
\end{equation}
where \(\pi_\Gamma\) is the restriction of \(\pi\) to \(\Gamma\).
Our standard Riemann-Hurwitz argument shows
\begin{equation}\label{prop1:e8}
2 g(\Gamma) - 2 + \left|\sigma(\Gamma)\right|\ge 
\gamma\left(2 g(G) - 2 + \deg B - \deg B_1 + c\right).
\end{equation}

If \(\supp \pi(\Gamma) = G''\), then
\begin{equation}\label{prop1:e9}
\pi_\Gamma^{-1}\left(G''\cap \left(G \cup \f\right)\right)
= \pi_\Gamma^{-1}\left(G''\cap \left(W_0^{(2)} \cup W_0^{(3)}\right)\right)
\subset \sigma(\Gamma)
\end{equation}
and by Riemann-Hurwitz, we have
\begin{equation}\label{prop1:e10}
\begin{split}
2 g(\Gamma) - 2 + \left|\sigma(\Gamma)\right|\ge 
\gamma(& 2 g(G) - 2\\
& + \deg B - \deg B_1 - \deg B_l + c).
\end{split}
\end{equation}

If $\pi(\Gamma)$ meets $W_0$ properly and dominates $G$
via \(f\), then by the
induction hypothesis \eqref{prop1:e5} and Riemann-Hurwitz,
\begin{equation}\label{prop1:e11}
\begin{split}
&\quad 2 g(\Gamma) - 2 + |\sigma(\Gamma)|\\
&\ge \gamma \left(2g(G) - 2 + \deg B - \deg B_1 - \max_{2\le k\le l-1}\deg B_k
+ c\right).
\end{split}
\end{equation}

In summary, the combination of \eqref{prop1:e8}, \eqref{prop1:e10} and
\eqref{prop1:e11} produces
\begin{equation}\label{prop1:e12}
\begin{split}
2 g(\Gamma) - 2 + \left|\sigma(\Gamma)\right|\ge 
\gamma\bigg(& 2 g(G) - 2 + \deg B - \deg B_1\\
&- \max_{2\le k\le l} \deg B_k + c\bigg)
\end{split}
\end{equation}
for each component $\Gamma\subset Y_0$ that dominates $G$ via the map
$f\circ \pi$ of degree $\gamma$. Therefore,
\begin{equation}\label{prop1:e13}
\begin{split}
\sum_{G\subset (f\circ \pi)(\Gamma)}
\left(2 g(\Gamma) - 2 + \left|\sigma(\Gamma)\right|\right)& \ge 
\beta\bigg( 2 g(G) - 2 + \deg B\\
&- \deg B_1 - \max_{2\le k\le l} \deg B_k + c\bigg)
\end{split}
\end{equation}
where we sum over all irreducible
components $\Gamma$ that dominate $G$ via the map $f\circ \pi$.
So it remains to show
\begin{equation}\label{prop1:e14}
\sum_{G\not\subset (f\circ \pi)(\Gamma)}
\left(2 p_a(\Gamma) - 2 + \left|\sigma(\Gamma)\right|\right) \ge 0
\end{equation}
where $\Gamma$ runs over all irreducible components of $Y_0$ that are
contractible under the map $f\circ \pi$, i.e., \(\pi(\Gamma)\) is
supported on a fiber of \(f: R\to G\).

The combination of
\eqref{prop1:e13} and \eqref{prop1:e14} will produce
\begin{equation}\label{prop1:e15}
\begin{split}
\sum_{\Gamma\subset Y_0}
\left(2 p_a(\Gamma) - 2 + \left|\sigma(\Gamma)\right|\right)& \ge 
\beta\bigg( 2 g(G) - 2 + \deg B\\
&- \deg B_1 - \max_{2\le k\le l} \deg B_k + c\bigg)
\end{split}
\end{equation}
and this implies \eqref{prop1:e6} thanks to \eqref{s2:ss2:e2}.

Let $Z$ be a maximal connected component of
\begin{equation}\label{prop1:e16}
\bigcup_{G\not\subset (f\circ \pi)(\Gamma)} \Gamma \subset Y_0.
\end{equation}
In order to prove \eqref{prop1:e14}, it suffices to prove
\begin{equation}\label{prop1:e17}
\sum_{\Gamma\subset Z}
\left(2 p_a(\Gamma) - 2 + \left|\sigma(\Gamma)\right|\right) \ge 0
\end{equation}
for all such $Z$.

We extend the notation $\sigma(\cdot)$ to \(Z\), for which \(\sigma(Z)\)
is the finite subset of $Z$ consisting of the marked points $Z\cap\lim_{t\to 0}
\left(Y_t\cap \pi^{-1}(W)\right)$ and the intersections between $Z$
and the rest of the components of $Y_0$. Then we have
\begin{equation}\label{prop1:e18}
2p_a(Z) - 2 + \left|\sigma(Z)\right| =
\sum_{\Gamma\subset Z}
\left(2 p_a(\Gamma) - 2 + \left|\sigma(\Gamma)\right|\right).
\end{equation}
Obviously, $p_a(Z)\ge 0$ since $Z$ is connected and $\sigma(Z)$
consists of at least one point unless $Z = Y_0$, which is impossible
since we assume that $Y_0$ dominates $G$ under the map $f\circ \pi$. So
\begin{equation}\label{prop1:e19}
\sum_{\Gamma\subset Z}
\left(2 p_a(\Gamma) - 2 + \left|\sigma(\Gamma)\right|\right)
= 2p_a(Z) - 2 + \left|\sigma(Z)\right|\ge -1.
\end{equation}
Suppose that the equality in \eqref{prop1:e19} holds. Then we must have
$p_a(Z) = 0$ and $\left|\sigma(Z)\right| = 1$. This means that every
irreducible component of $Z$ is smooth and rational, the dual graph of
$Z$ is a tree, $Z$ meets the rest of $Y_0$ at exactly one point and
there are no marked points on $Z$. In addition, $Z$ cannot be
contracted to a point by $\pi$ due to our minimality assumption on
$\pi: Y\to X$ and hence $\pi(Z)$ is supported along a fiber of the
projection $f: R\to G$.

For each component $W^{(k)}\subset W$ ($1\le k\le 3$), we write
\begin{equation}\label{prop1:e20}
\pi^*\left(W^{(k)}\right) = L_k + M_k + N_k
\end{equation}
where $L_k$ consists of the sections of $Y\to\Delta$,
$M_k$ is supported along $Z$ and $N_k$ is supported along the rest of
$Y_0$. Since there are no marked points on $Z$, $L_k\cdot Z = 0$ 
and $L_k\cdot M_k = 0$ for $k=1,2,3$.

Let $\Gamma$ be the only irreducible component of $Y_0$ that
intersects $Z$. If $\pi(\Gamma)$ is supported on $G$,
we have $\Gamma\not\subset N_1$ and set $k = 1$; otherwise, if
$\pi(\Gamma)$ is not supported on $G$, we have $\Gamma\not\subset N_3$
and set $k = 3$. In any event, we have $N_k\cdot Z = 0$, $N_k\cdot M_k =
0$ and $M_k^2 = M_k(L_k + M_k + N_k)\ge 0$. This is possible only if $M_k
= \emptyset$. But then we have $Z(L_k + M_k + N_k) = 0$, which
contradicts the fact that $\pi(Z)$ is supported along a fiber of $f:
R\to G$ and hence 
\(Z\cdot\pi^*\left(W^{(k)}\right) = \pi(Z)\cdot W^{(k)} > 0\).

Hence we must have $2p_a(Z) - 2 + \left|\sigma(Z)\right|\ge 0$ and
\eqref{prop1:e17} follows. This finishes the proof of the proposition.
\end{proof}

\section{Algebraic Hyperbolicity of Log Surfaces}\label{s3}

In this section, we will show some applications of \thmref{t3} by proving
\coref{cor0}-\ref{cor2} and \thmref{t4}.

\subsection{Proof of \coref{cor0}-\ref{cor2}}\label{s3:ss1}

\begin{proof}[Proof of \coref{cor0}]
The first inequality \eqref{cor0:e1} follows immediately from
\thmref{t3} and the second \eqref{cor0:e2} follows from
\eqref{cor0:e1} and Mori's cone theorem, as mentioned before. It
remains to verify the last statement, which requires us to show that
\(K_S + 2B\) is NEF if \((K_S + 2B) B\ge 0\).

By Mori's cone theorem, it suffices to show that \((K_S + 2B) G \ge
0\) for all rational curves \(G\subset S\)
with \(-3\le K_S\cdot G < 0\). Note that if such \(G\) exists, \(S\)
is rational and covered by the curves in the numerical class of
\(G\). So it is enough to show that \(K_S + 2B\) is effective.
To see this, let \(f: \wt{S}\to S\) be a minimal desingularization of \(S\).
Then \(K_{\wt{S}} = f^* K_S\) since \(S\) is canonical. So \(K_S +
2B\) is effective if and only if \(K_{\wt{S}} + 2\wt{B}\) is, where
\(\wt{B} = f^* B\). So it comes down to a Riemann-Roch computation on
\(\wt{S}\):
\begin{equation}\label{cor0:e4}
h^0(\CO_{\wt{S}}(K_{\wt{S}} + 2\wt{B})) \ge (K_{\wt{S}} +
2\wt{B})\wt{B} + 1 = (K_S + 2B)B + 1 \ge 1
\end{equation}
and we are done.
\end{proof}

Both \coref{cor1} and \ref{cor2} are direct consequences of
\thmref{t3}.

\begin{proof}[Proof of \coref{cor1}]
Every BPF divisor on $\F_n$ lies in the cone of numerically effective
divisors (NEF cone) generated by $F$ and
$M + nF$ (BPF is equivalent to NEF on $\F_n$). Hence each $D_k$ is
linearly equivalent to a sum of divisors $F$ and $M + nF$.
Applying \thmref{t3}, we may take $\varphi(C) = \deg C$, each
$B_i$ to be either $F$ or $M+nF$ and $\epsilon$ to be the smallest
of the following numbers:
\begin{equation}\label{cor1:e2}
\begin{split}
\min_C \frac{D\cdot C - 2}{\deg C} &= \min(a - 2, b - an - 2)\\
\min_C \frac{(K_{\F_n} + D - F) C}{\deg C} &= \min(a - 2, b - an +
n - 3)\\
\min_C \frac{(K_{\F_n} + D - M - nF) C}{\deg C}
&= \min(a - 3, b - an + n - 2)
\end{split}
\end{equation}
where $C$ runs over all effective divisors of $\F_n$; it is enough to
check \eqref{cor1:e2}
with $C\in \{M, F\}$ since the cone of effective divisors of $\F_n$ is
generated by $M$ and $F$. And \eqref{cor1:e1} follows immediately.
\end{proof}

\begin{proof}[Proof of \coref{cor2}]
Again BPF is equivalent to NEF on $\wt{\P}^2$.
The NEF cone of $\wt{\P}^2$ is generated by divisors $M$ with the
property that $M\cdot R_i = 0$ or $1$ for each $R_i$. 
Applying \thmref{t3}, we may take $\varphi(C) = \deg C$ and $\epsilon$
to be the smaller of the following two numbers:
\begin{equation}\label{cor2:e2}
\min_C \frac{D\cdot C - 2}{\deg C} \text{ and }
\min_{C, M} \frac{(K_{\wt{\P}^2} + D - M) C}{\deg C}
\end{equation}
where $M$ runs over all divisors satisfying $M\cdot R_i = 0$ or $1$
for each $R_i$ and $C$ runs over all effective divisors on
$\wt{\P}^2$. It is enough to check \eqref{cor2:e2} with $C\in \{R_1,
R_2, ..., R_n\}$ since the cone of effective divisors on $\wt{\P}^2$ is
generated by $R_1, R_2, ..., R_n$ and then it is obvious that
$\epsilon = \min D\cdot R_i - 2$.
\end{proof}

\subsection{Complements of curves on K3 surfaces}\label{s3:ss2}

Next, we will prove \thmref{t4}.

We have to degenerate $S$ and $D$ simultaneously as
mentioned at the beginning. Let us consider K3 surfaces with Picard lattice
\begin{equation}\label{s3:ss2:e0}
\begin{pmatrix}
0 & 2\\
2 & 0
\end{pmatrix}
\end{equation}
or
\begin{equation}\label{s3:ss2:e1}
\begin{pmatrix}
0 & 2\\
2 & -2
\end{pmatrix}.
\end{equation}
Such K3 surfaces can be realized as double covers of a rational surface $R$
ramified along a curve in $\P H^0(\CO_R(-2K_R))$,
where $R = \F_0 = \P^1\times\P^1$
for K3 surfaces with Picard lattice \eqref{s3:ss2:e0} and $R = \F_1$
for K3 surfaces with Picard lattice \eqref{s3:ss2:e1}. Let $M$ and $F$
be the generators of the Picard group of such a K3 surface, where
$M\cdot F = 2$, $F^2 = 0$ and $M^2 = 0$ or $-2$.

Let $X$ be a family of surfaces over the disk $\Delta$
parameterized by $t$, where $X_t$ is a general K3 surface of
genus $g$ for $t\ne 0$ and the central fiber $X_0$ is a K3 surface
with Picard lattice \eqref{s3:ss2:e0} or \eqref{s3:ss2:e1}. The limit
of the primitive divisors of $X_t$ is $M + \lfloor g/2\rfloor F$ on
$X_0$ and the Picard lattice of $X_0$ is \eqref{s3:ss2:e0} if $g$ is odd and
\eqref{s3:ss2:e1} if $g$ is even.

Let $W$ be a divisor on $X$ such that $W_t$ is a curve
in the primitive class of $X_t$ for $t\ne 0$ and $W_0 = G$ is a general
member in the linear series $\P H^0(\CO_{X_0}(M + \lfloor g/2\rfloor F))$.

Let $Y$ be a flat family of curves over $\Delta$ with the commutative
diagram \eqref{t5:e1}. As before, we assume that $Y$ is smooth and
irreducible, the map $Y_t\to X$ is generically 1-1 onto its image for
$t\ne 0$,
$\pi(Y)$ meets $W$ properly in $X$ and $Y\to X$ is a family of
semistable maps with marked points $Y_t\cap \pi^{-1}(W)$ on $Y_t$
($t\ne 0$).
Our goal is to prove that
\begin{equation}\label{s3:ss2:e2}
\begin{split}
2 g(\pi(Y_t)) - 2 + i_{X}(\pi(Y_t), W) &\ge
\frac{\lfloor (g-3)/2\rfloor}{g-1}(W \cdot \pi(Y_0))\\
&= \lfloor (g-3)/2\rfloor (F\cdot \pi(Y_0))
\end{split}
\end{equation}
for $t\ne 0$, where the intersection \(F\cdot \pi(Y_0)\) is taken on the
surface $X_0$.

Let $\sigma(\Gamma)\subset\Gamma$ be defined as in the proof of
\thmref{t3} for each component $\Gamma\subset Y_0$ and we have
\begin{equation}\label{s3:ss2:e3}
2 g(\pi(Y_t)) - 2 + i_{X}(\pi(Y_t), W) = \sum_{\Gamma\subset Y_0}
\left(2p_a(\Gamma) - 2 + |\sigma(\Gamma)|\right).
\end{equation}
Therefore, we can reduce \eqref{s3:ss2:e2} to showing that
\begin{equation}\label{s3:ss2:e4}
2p_a(\Gamma) - 2 + |\sigma(\Gamma)|
\ge \lfloor (g-3)/2\rfloor (F\cdot \pi(\Gamma))
\end{equation}
for each component $\Gamma\subset Y_0$.

The inequality \eqref{s3:ss2:e4}
obviously holds for $\Gamma$ contractible under
$\pi$. So we assume that $\Gamma\subset Y_0$ is noncontractible.

If $\Gamma$ dominates $W_0 = G$, then
\begin{equation}\label{s3:ss2:e5}
2 g(\Gamma) - 2 \ge \gamma(2g(G) - 2) = (g - 1) (F\cdot \pi(\Gamma))
\end{equation}
by Riemann-Hurwitz, where $\gamma$ is the degree of the map $\Gamma\to
G$. And hence \eqref{s3:ss2:e4} holds for $\Gamma$ dominating $G$.

Suppose that $\Gamma$ is a component of $Y_0$ such that $\pi(\Gamma)$
meets $G$ properly. It is not hard to see that
\begin{equation}\label{s3:ss2:e6}
\pi_\Gamma^{-1}\left(C\cap G\right)\subset \sigma(\Gamma)
\end{equation}
where $C = \left(\pi(\Gamma)\right)_\red$ is the reduced image of
\(\Gamma\) under \(\pi\) and
\(\pi_\Gamma: \Gamma\to C\) is the restriction of \(\pi\) to
\(\Gamma\).
By Riemann-Hurwitz, we have
\begin{equation}\label{s3:ss2:e7}
2 g(\Gamma) - 2 + \left|\sigma(\Gamma)\right| \ge \gamma\left(2g(C) -
2 + i_{X_0}(C, G)\right)
\end{equation}
where $\gamma$ is the degree of the map $\Gamma\to C$.
Hence in order to show \eqref{s3:ss2:e4}, it is enough to show
\begin{equation}\label{s3:ss2:e8}
2 g(C) - 2 + i_{X_0}(C, G) \ge \lfloor (g-3)/2\rfloor (F\cdot C)
\end{equation}
for each reduced irreducible curve $C\subset X_0$ with $C\ne G$.
This naturally leads us to apply \thmref{t3} to $(X_0, G)$ with
the following argument.

We take $\varphi(C) = F\cdot C$ in \thmref{t3}.

If $g$ is odd, we break up $G$ into a union of one curve in $\P
H^0(\CO(M))$ and $\lfloor g/2\rfloor$ curves in $\P H^0(\CO(F))$,
i.e., we take $B_i = M$ or
$F$ in \thmref{t3}. It is easy to check that
\begin{equation}\label{s3:ss2:e9}
(K_{X_0} + G - M)C = \lfloor g/2\rfloor (F\cdot C),
\end{equation}
\begin{equation}\label{s3:ss2:e10}
(K_{X_0} + G - F)C \ge \left(\lfloor g/2\rfloor - 1\right) 
(F\cdot C)
\end{equation}
and
\begin{equation}\label{s3:ss2:e11}
G\cdot C - 2\ge \left(\lfloor g/2\rfloor - 1\right) (F\cdot C)
\end{equation}
for all curves $C\subset X_0$. And hence \eqref{s3:ss2:e8} follows.

If $g$ is even, we break up $G$ into a union of exactly one curve in
$\P H^0(\CO(M + F))$
and $\lfloor g/2\rfloor - 1$ curves in $\P H^0(\CO(F))$, i.e.,
we take $B_i = M +
F$ or $F$ in \thmref{t3}. It is easy to check that
\begin{equation}\label{s3:ss2:e12}
(K_{X_0} + G - M - F) C = (\lfloor g/2\rfloor - 1) (F\cdot C),
\end{equation}
\begin{equation}\label{s3:ss2:e13}
(K_{X_0} + G - F) C \ge \left(\lfloor g/2\rfloor - 2\right) 
(F\cdot C)
\end{equation}
and
\begin{equation}\label{s3:ss2:e14}
G\cdot C - 2\ge \left(\lfloor g/2\rfloor - 2\right) (F\cdot C)
\end{equation}
for all curves $C\subset X_0$. And hence \eqref{s3:ss2:e8}
follows. This finishes the proof of \eqref{s3:ss2:e4} and consequently
the proof of \thmref{t4} for $P = \emptyset$. Little in the above
proof needs to be changed for the case $P\ne \emptyset$. We will leave
that to the readers.

\section{Algebraic Hyperbolicity of Projective Surfaces}\label{s4}

\subsection{Proof of \thmref{t6}}\label{s4:ss1}

First, we will show how to derive \thmref{t6} from \thmref{t2} and
\ref{t5}, while the proof of the latter will be postponed to
\ssecref{s4:ss3}.

Obviously, what we are supposed to do is to degenerate each \(V_j\) to
a union of \(\sum_{i=1}^m a_{ij}\) hypersurfaces consisting of 
\(a_{ij}\) hypersurfaces in \(\P \L_i\) for \(i = 1,2,...,m\).

Let \(Z^{(j)}\subset W\times\Delta\) be a pencil of hypersurfaces in $W$
whose general fiber \(Z_t^{(j)}\) are general members of the linear
system \eqref{t6:e1} and
whose central fiber \(Z_0^{(j)}\) is the union of \(\sum_{i=1}^m a_{ij}\)
hypersurfaces consisting of \(a_{1j}\) hypersurfaces in \(\P \L_1\),
\(a_{2j}\) hypersurfaces in \(\P \L_2\), ..., \(a_{mj}\) hypersurfaces
in \(\P \L_m\). Let \(X = \cap_{j=1}^{n} Z^{(j)}\subset
W\times\Delta\). Obviously, \(X_0\) is of normal crossing and \(X\)
has only rational double points at \(Q = X_{sing}\) for a general
choice of \(Z^{(j)}\). Also since \(W\) is smooth in codimension 2,
\(X_0\) can be chosen away from \(W_{sing}\). So we may simply assume
\(W\) to be smooth.

Let \(Y\) be a flat family of curves over
\(\Delta\) with the commutative diagram \eqref{t5:e1}. We want to show
that \eqref{t6:e4} holds for \(C = \pi(Y_t)\) (assuming that \(\pi:
Y_t \to X\) is birational onto its image for \(t\ne 0\)). Then by
\thmref{t5}, it is enough to show
\begin{equation}\label{t6:e6}
2 g_Q^{vir}(\pi(Y_0)) - 2 \ge \epsilon \varphi(\pi(Y_0)).
\end{equation}
Let \(\Gamma\) be an irreducible component of \(\pi(Y_0)\).
Our goal is to show
\begin{equation}\label{t6:e7}
\phi_Q(\Gamma) \ge \epsilon \varphi(\Gamma)
\end{equation}
which should imply \eqref{t6:e6} immediately by the definition of
virtual genus, where \(\phi_Q(\Gamma)\) is the notation used in
\defnref{defn4}.

Let us first check the case \(\Gamma\subset \partial D\) for some
component \(D\subset X_0\). In this case, \(\Gamma = U_1\cap U_2\cap
... \cap U_n\cap U_{n+1}\), where each \(U_j\) is a member of \(\P
\L_{i}\) for some \(1\le i\le m\), \(U_j\subset Z_0^{(j)}\)
for \(j\le n\) and \(U_{n+1}\subset Z_0^{(l)}\) for some \(l\).
Then
\begin{equation}\label{t6:e8}
\begin{split}
\phi_Q(\Gamma) &= 2 g(\Gamma) - 2 + \left(\sum_{j=1}^n Z_0^{(j)} -
\sum_{k=1}^{n+1} U_k\right)\Gamma\\
&= \left(K_W + \sum_{k=1}^{n+1} U_k\right)\Gamma + \left(\sum_{j=1}^n
Z_0^{(j)} - \sum_{k=1}^{n+1} U_k\right)\Gamma\\
&= \left(K_W + \sum_{i, j} a_{ij} B_i\right)\Gamma \ge
\epsilon \varphi(\Gamma).
\end{split}
\end{equation}
Note that by B\'ezout, \(\Gamma = U_1\cap U_2\cap
... \cap U_n\cap U_{n+1}\) is smooth. However, \(\Gamma\) could very well
be disconnected. But the monodromy action on the connected components
of \(\Gamma\) as \(U_j\) varies is transitive.
In particular, this means that the components of \(\Gamma\) are
deformationally equivalent to each other.
Therefore,
\eqref{t6:e8} holds for each irreducible component of \(\Gamma\),
which is what we really need.

Now let us deal with the case \(\Gamma\subset D\) for some component
\(D\subset X_0\) and \(\Gamma\not\subset \partial D\). Let
\(D = U_1\cap U_2\cap
... \cap U_n\), where each \(U_j\) is a member of \(\P
\L_i\) for some \(1\le i\le m\) and \(U_j\subset Z_0^{(j)}\).
As expected, we want to apply \thmref{t2} to \((D, \partial D,
Q\cap \partial D)\). There are a few things to verify.
First, we need to check
\begin{equation}\label{t6:e9}
(K_D + \partial D - G)\Gamma\ge \epsilon\varphi(\Gamma)
\end{equation}
where \(G = D\cap U\subset \partial D\) is cut out by \(U\subset
Z_0^{(j)}\) and \(U\in \P \L_i\) for some \(i\) and \(j\).
Observe that \(K_D
+ \partial D = K_W + \sum_{i, j} a_{ij} B_i\) and \eqref{t6:e9}
follows immediately from \eqref{t6:e2}.
Second, we need to check
\begin{equation}\label{t6:e10}
2 g(\Gamma) - 2 + \partial D \cdot \Gamma \ge \epsilon\varphi(\Gamma).
\end{equation}
This follows from \eqref{t6:e3} since \(\partial D \cdot \Gamma
= \left(\sum_{i, j} a_{ij} B_i - \sum_{k=1}^n U_k\right)
\Gamma\). Finally, let \(Q\cap \partial D \subset F\cap \partial D\)
for some curve \(F\subset D\) and we need to show
\begin{equation}\label{t6:e11}
2 g(\Gamma) - 2 + i_D(\Gamma, \partial D, Q\cap \partial D) \ge
\epsilon\varphi(\Gamma)
\end{equation}
for each irreducible component \(\Gamma\subset F\).
Actually, \(Q\) and \(F\) can
be described as follows: \(F = F_1 \cup F_2 \cup ...\cup F_n\) and
\begin{equation}\label{t6:e12}
Q\cap \partial D = \bigcup_{j=1}^n 
\left(F_j\cap (Z_0^{(j)}\backslash U_j)\right)
\end{equation}
where \(F_j\) is cut out on \(D\) by a general member \(Z_t^{(j)}\)
of the pencil \(Z^{(j)}\), i.e., it is cut out by a general member in
the linear series \eqref{t6:e1}. We have
\begin{equation}\label{t6:e13}
\begin{split}
&\quad 2 g(F_j) - 2 + i_D(F_j, \partial D, Q\cap \partial D)\\
&= (K_D + Z_t^{(j)}) F_j + (\partial D - Z_0^{(j)} + U_j) F_j\\
&= \left(K_W + \sum_{i, j} a_{ij} B_i + U_j\right) F_j \ge
\epsilon\varphi(F_j).
\end{split}
\end{equation}
Again, \(F_j\) may very well be disconnected but the monodromy action
on its connected components is transitive. Therefore, \eqref{t6:e11}
holds for every irreducible component \(\Gamma\subset F_j\) by
\eqref{t6:e13}. 

Combining \eqref{t6:e9}-\eqref{t6:e11} and applying \thmref{t2} to
\((D, \partial D, Q\cap \partial D)\), we conclude
\begin{equation}\label{t6:e14}
\phi_Q(\Gamma) = 2 g(\Gamma) - 2 + i_D(\Gamma, \partial D, Q\cap
\partial D) \ge \epsilon\varphi(\Gamma).
\end{equation}
This finishes the proof of \eqref{t6:e7} and hence \thmref{t6}.

\subsection{Proof of \thmref{t7}}\label{s4:ss2}

Next, we will show how to derive \thmref{t7} from \thmref{t4} and
\ref{t5}.

Suppose that $S$ is a complete intersection of type $(2, 2, 3)$ in
$\P^5$. We can split $S$ into a union of a complete intersection of type
$(2, 2, 1)$ and a complete intersection of type $(2,2,2)$.

Let $X\subset\P^5\times\Delta$ be a family of complete intersections
of type $(2,2,3)$ with central fiber $X_0 = R\cup W$, where $R$ and $W$
are complete intersections of type $(2,2,1)$ and $(2,2,2)$,
respectively. That is, $R$ is a Del Pezzo surface of degree $4$, i.e.,
the blowup of $\P^2$ at five points in general position and $W$ is a
K3 surface of genus $5$. Let $E = R\cap W$ and $X_{sing} = Q = \{q_1,
q_2, ..., q_{24}\}\subset E$ be the set of rational double points of
$X$. Obviously, $E$ is in the primitive class of $W$ and
$E\in \P H^0\left(\CO_R(-2K_R)\right)$
on $R$; $Q$ is cut out on $E$ by a curve in
$\P H^0\left(\CO_W(3E)\right)$ on $W$ and by a
curve in $\P H^0\left(\CO_R(-3K_R)\right)$ on $R$.

Let $Y$ be a flat family of curves with the commutative diagram
\eqref{t5:e1}. Our goal is to prove that
\begin{equation}\label{s4:ss2:e0}
2 g(Y_t) - 2\ge \frac{1}{6}\deg \pi(Y_t)
\end{equation}
for $t\ne 0$. Again by \thmref{t5}, it suffices to show that
\begin{equation}\label{s4:ss2:e1}
2 g_Q^{vir}(\pi(Y_0)) - 2 \ge \frac{1}{6} \deg \pi(Y_0).
\end{equation}
Let $\pi(Y_0) = \Sigma_R \cup \Sigma_W\cup mE$, where $\Sigma_R\subset
R$, $\Sigma_W\subset W$, $E\not\subset \Sigma_R, \Sigma_W$ and $m =
\mu(E)$ is the multiplicity of $E$ in $\pi(Y_0)$. By
\thmref{t4}, $(W, E, Q)$ is algebraically hyperbolic and
\begin{equation}\label{s4:ss2:e2}
2 g(\Gamma) - 2 + i_W(\Gamma, E, Q)\ge \frac{1}{4} \deg \Gamma
\end{equation}
for each irreducible component
$\Gamma\subset \Sigma_W$. We are done if $(R, E, Q)$ is
algebraically hyperbolic too. However, it is not since obviously every $-1$
rational curve on $R$ meets $E$ at two points. On the other hand, we do have
\begin{equation}\label{s4:ss2:e3}
2 g(\Gamma) - 2 + i_R(\Gamma, E, Q)\ge 0
\end{equation}
for each irreducible
component $\Gamma\subset\Sigma_R$ by \coref{cor2}. Note that
\begin{equation}\label{s4:ss2:e4}
\begin{split}
2 g_Q^{vir}(\pi(Y_0)) &- 2 = \sum_{\Gamma\subset \Sigma_R}
\mu(\Gamma)\left(2 g(\Gamma) - 2 + i_R(\Gamma, E, Q)\right)\\
&+ \sum_{\Gamma\subset \Sigma_W}
\mu(\Gamma)\left(2 g(\Gamma) - 2 + i_W(\Gamma, E, Q)\right) + m \deg E
\end{split}
\end{equation}
where $\mu(\Gamma)$ is the multiplicity of a component $\Gamma$ in
$\pi(Y_0)$. Hence \eqref{s4:ss2:e1} will hold unless $\deg \Sigma_R$ is
large compared with $\deg\Sigma_W$ and $m$. To make this precise, let
$\deg \Sigma_R = r$, $\deg \Sigma_W = w$ and
\begin{equation}\label{s4:ss2:e5}
\beta = \left\lceil \frac{r}{12} - \frac{w}{24} - \frac{m}{3}
\right\rceil.
\end{equation}
A standard analysis of limit linear series on $X_0$ shows that if
$\beta > 0$, $\Sigma_R$ has a singularity of multiplicity at least
$\beta$ at each $q_i\in Q$.

Since $Q$ is cut out on $E$ by a general curve in the linear system
$\P H^0(\CO_R(-3K_R))$,
$\sum_{i=1}^{24} q_i$ moves in a complete linear series of degree $24$
on $E$. And
since $g(E) = 5$, there exists a pencil in the linear series
$\P H^0\left(\CO_E\left(\sum_{i=1}^{24}
q_i\right)\right)$ with $18$ base points, which means that we may fix $16$ points
among $Q$ and let the rest vary. This is where we may apply the
argument in \remref{rem6}.
Suppose that $\Gamma\subset\Sigma_R$ is a component of
$\Sigma_R$ passing through $l$ points among $Q$, say that $\Gamma$
passes through $q_1, q_2, ..., q_l\in Q$. Let $m_1, m_2, ..., m_l$ be the
multiplicities of $\Gamma$ at $q_1, q_2, ..., q_l$. Without the loss
of generality, we assume that $m_1\ge m_2\ge ...\ge m_l$.
If $l > 18$, we fix $\{q_1, q_2, ..., q_{18}\}$ and vary $\{q_{19},
q_{20}, ..., q_{l}\}$ and then the argument in \remref{rem6} yields
\begin{equation}\label{s4:ss2:e6}
\begin{split}
2 g(\Gamma) - 2 + i_R(\Gamma, E, Q) &\ge -K_R\cdot\Gamma -
\sum_{i=19}^l m_i\\
&\ge -K_R \cdot \Gamma - \frac{l-18}{l} E\cdot \Gamma\\
&=\frac{36 - l}{l} \deg\Gamma \ge \frac{1}{2} \deg\Gamma
\end{split}
\end{equation}
where we observe
\begin{equation}\label{s4:ss2:e7}
\sum_{i=19}^l m_i \le \frac{l-18}{l} E\cdot \Gamma
\end{equation}
since $\sum_{i=1}^l m_i \le E\cdot \Gamma$ and $m_1\ge m_2\ge ...\ge m_l$.
If $1\le l \le 18$, we fix $\{q_1, q_2, ..., q_{l-1}\}$ and vary $q_l$
and then the same argument in \remref{rem6} yields
\begin{equation}\label{s4:ss2:e8}
\begin{split}
2 g(\Gamma) - 2 + i_R(\Gamma, E, Q) &\ge -K_R\cdot\Gamma - m_l\\
&\ge -K_R \cdot \Gamma + \frac{1}{2} K_R\cdot \Gamma =
\frac{1}{2}\deg\Gamma
\end{split}
\end{equation}
where we observe that $-K_R\cdot \Gamma \ge 2 m_l$ since
the Seshadri constant of $-K_R$ on $R$ is $2$.
Combining \eqref{s4:ss2:e6} and \eqref{s4:ss2:e8}, we see that
\begin{equation}\label{s4:ss2:e9}
2 g(\Gamma) - 2 + i_R(\Gamma, E, Q) \ge \frac{1}{2} \deg\Gamma
\end{equation}
if $\Gamma\subset\Sigma_R$ passes through at least one point of
$Q$.

Let $Z\subset \Sigma_R$ be the union of the components that pass
through at least one point of $Q$ and let $\deg Z = z$. Obviously,
\begin{equation}\label{s4:ss2:e10}
r\ge z \ge 12\beta\ge r - \frac{w}{2} - 4m.
\end{equation}
Combining \eqref{s4:ss2:e2}, \eqref{s4:ss2:e3}, \eqref{s4:ss2:e4} and
\eqref{s4:ss2:e9}, we have
\begin{equation}\label{s4:ss2:e11}
2 g_Q^{vir}(\pi(Y_0)) - 2 \ge \frac{w}{4} + \frac{z}{2} + 8m.
\end{equation}
Obviously,
\begin{equation}\label{s4:ss2:e12}
\deg \pi(Y_0) = r + w + 8m.
\end{equation}
It remains to minimize the RHS of \eqref{s4:ss2:e11} under the constraints
\eqref{s4:ss2:e10}, \eqref{s4:ss2:e12} and $r, w, z, m\ge 0$ (with $\deg
\pi(Y_0)$ fixed) and \eqref{s4:ss2:e1} follows. This finishes the proof of
\thmref{t7} in the case that $S\subset\P^5$ is a complete intersection
of type $(2, 2, 3)$.

For a complete intersection $S\subset\P^6$ of type $(2, 2, 2, 2)$,
we can degenerate it to a union of two surfaces both of type $(2, 2, 2,
1)$. Let $X\subset\P^6\times\Delta$ be the corresponding family of
surfaces and $X_0 = W_1 \cup W_2$, where $W_k\subset \P^6$ is a
complete intersection of type $(2, 2, 2, 1)$ for $k = 1, 2$.
Let $E = W_1\cap W_2$ and $X_{sing} = Q$. Obviously, $W_k$ is a K3
surface of genus $5$, $E\subset W_k$ is a curve in the primitive
class of $W_k$ and $Q$ is cut out on $E$ by a curve in
$\P H^0\left(\CO_{W_k}(2E)\right)$ for $k= 1,2$. 

Let $Y$ be a flat family of curves with the commutative diagram
\eqref{t5:e1}. Again it suffices to show
\begin{equation}\label{s4:ss2:e12.1}
2 g_Q^{vir}(\pi(Y_0)) - 2 \ge \frac{1}{4} \deg \pi(Y_0).
\end{equation}
By \thmref{t4}, $(W_k, E, Q)$ is algebraically hyperbolic and
\begin{equation}\label{s4:ss2:e13}
2 g(\Gamma) - 2 + i_{W_k}(\Gamma, E, Q)\ge \frac{1}{4} \deg \Gamma
\end{equation}
for each component $\Gamma$ of $\pi(Y_0)$ with $\Gamma\subset W_k$ and
$\Gamma\ne E$. And it is obvious that
\begin{equation}\label{s4:ss2:e14}
2 g(\Gamma) - 2 \ge \deg\Gamma
\end{equation}
for each component $\Gamma\subset \pi(Y_0)$ with $\Gamma = E$. Then
\eqref{s4:ss2:e12.1} follows immediately from \eqref{s4:ss2:e13} and
\eqref{s4:ss2:e14}.

\subsection{Proof of \thmref{t5}}\label{s4:ss3}

\lemref{lem2} deals with the local behavior of \(X\) and \(Y\) at a
point \(p\in \partial D_j\)
where \(X\) is smooth at \(p\). The following lemma
will tell us something about \(X\) and \(Y\) at a point \(q\in Q\)
where \(X\) is singular.

\begin{lem}\label{lem3}
Let \(X\subset \Delta_{x_1 x_2 ... x_r}^{r}\times \Delta\) be the
hypersurface given by \(x_1x_2...x_n = t f(t, x_1, x_2, ..., x_r)\),
where \(n < r\), \(f(0, 0, 0, ..., 0) = 0\) and
\(f(t, x_1, x_2, ..., x_r) \not\equiv 0\) along
\(x_1 = x_2 = ... = x_n = t = 0\). Let \(X_0 = D = \cup_{k=1}^n D_k\)
where \(D_k = \{ x_k = t = 0 \}\). And let \(Q\) be the singular locus
of \(X\), i.e., \(Q\) is cut out on \(X\) by \(f(t, x_1, x_2, ...,
x_r) = 0\).

Let
$Y$ be a flat family of curves over $\Delta$ with the commutative diagram
\eqref{t5:e1}. Suppose that there exists \(J\subset \{1,2,...,n\}\)
and a reduced irreducible curve \(\Gamma\subset \pi(Y_0)\)
such that \(\Gamma\subset
D_J\), \(\Gamma\not\subset\partial D_J\) and \(i_{D_J}(\Gamma,
\partial D_J, Q \cap \partial D_J) > 0\). Then there exists a
curve \(\Gamma'\subset \pi(Y_0)\) such that \(\Gamma' \subset
\cup_{k\not\in J} D_k\).
\end{lem}

\begin{proof}
Let $\wt{X}$ be the blowup of $X$ along $\cup_{j\in J} D_j$. The map
$f: \wt{X}\to X$ is a small morphism ($f_* \CO_{\wt{X}} = \CO_X$) and
it is an isomorphism away from
\begin{equation}\label{s4:ss3:e3}
Q\cap \left(\bigcup_{j\in J} D_j\right)\cap
\left(\bigcup_{k\not\in J} D_k\right).
\end{equation}
The central fiber of $\wt{X}$ looks like
\begin{equation}\label{s4:ss3:e4}
\wt{X}_0 = \left(\bigcup_{j\in J} \wt{D}_j\right)\cup
\left(\bigcup_{k\not\in J} D_k\right)
\end{equation}
where $f(\wt{D}_j) = D_j$ for $j\in J$. Let $\wt{D}_J = \cap_{j\in J}
\wt{D}_j$. It is not hard to see that $f: \wt{D}_J\to D_J$ is the
blowup of $D_J$ along $Q\cap \partial D_J$. We let 
\begin{equation}\label{s4:ss3:e5}
\partial \wt{D}_J
= \wt{D}_J\cap \left(\bigcup_{k\not\in J} D_k\right).
\end{equation}
Obviously, $\partial \wt{D}_J$ is the proper transform of 
\(\partial D_J\) under the blowup \(f: \wt{D}_J\to D_J\).

The rational map $Y\to \wt{X}$ can resolved into a regular map
$\wt{\pi}: \wt{Y}\to\wt{X}$ and we have the commutative diagram
\eqref{s2:ss3:e6}. Let \(\wt{\Gamma}\) be the proper transform of
\(\Gamma\) under the blowup \(f: \wt{D}_J\to D_J\).
Obviously, \(\wt{\Gamma}\subset \wt{\pi}(\wt{Y_0})\). Since \(i_{D_J}(\Gamma,
\partial D_J, Q \cap \partial D_J) > 0\), \(\wt{\Gamma}\) passes
through the close point \(q\) of \(\partial \wt{D}_J\). Apply
\lemref{lem2} to \((\wt{X}, \wt{Y}, \wt{\Gamma})\) locally at \(q\)
and we are done.
\end{proof}

Now let us go back to the proof of \thmref{t5}.

We may assume that $Y$ is smooth and irreducible and
$Y\to X$ is a family of semistable maps.
For each component $\Gamma\subset Y_0$, we define
$\sigma(\Gamma)\subset \Gamma$ to be the finite subset of \(\Gamma\)
consisting of the intersections between \(\Gamma\) and the
components of \(Y_0\) other than \(\Gamma\). Then
\begin{equation}\label{s4:ss3:e1}
2 g(Y_t) - 2 + i_X(\pi(Y_t)) = \sum_{\Gamma\subset
Y_0} \left( 2 p_a(\Gamma) - 2 + \left|\sigma(\Gamma)\right|\right).
\end{equation}

Let $\Gamma\subset Y_0$ be an irreducible component of $Y_0$ and
let $C = \left(\pi(\Gamma)\right)_\red$ be the reduced image of the
map $\Gamma\to X_0$. Our goal is to prove
\begin{equation}\label{s4:ss3:e2}
2 p_a(\Gamma) - 2 + \left|\sigma(\Gamma)\right| \ge \gamma \phi_Q(C)
\end{equation}
for each $\Gamma$, 
where $\gamma$ is the degree of the map $\pi: \Gamma\to C$. Here we
take \(\phi_Q(C) = 0\) if \(\Gamma\) is contractible under \(\pi\).
And \eqref{t5:e2} should follow easily from \eqref{s4:ss3:e1} and
\eqref{s4:ss3:e2}. Since \eqref{s4:ss3:e2} is trivial for
\(\Gamma\) contractible, let us assume that \(\Gamma\) is noncontractible.

Let $J\subset\{1,2,...,n\}$ be the index set such that
$C\subset D_j$ for $j\in J$ and $C\not\subset D_j$ for
$j\not\in J$. Let \(\wt{D}_J\) be the blowup of \(D\) along \(Q\cap
\partial D_J\) and \(\partial \wt{D}_J\) be the proper transform of
\(D_J\) under the blowup.

Let $\nu_\gamma: \Gamma^\nu\to\Gamma$ and
$\nu_c: C^\nu \to C$ be the normalizations of $\Gamma$
and $C$, respectively. We have the following commutative diagram
\begin{equation}\label{lem3:e1}
\begin{array}{ccc}
\Gamma^\nu & \mapright{\pi_\Gamma} & C^\nu\\
\mapdown{\nu_\gamma} & & \mapdown{\nu_c}\\
\Gamma & \mapright{\pi_\Gamma} & C
\end{array}
\end{equation}
where we use $\pi_\Gamma$ to denote both the maps
\(\Gamma^\nu\to C^\nu\) and \(\Gamma\to C\). Let \(\sigma(\Gamma^\nu)
= \nu_\gamma^{-1} (\sigma(\Gamma))\) and \(\sigma(C^\nu) =
f^{-1}(\partial \wt{D}_J)\), where \(f:
C^\nu\to \wt{D}_J\) is the map induced by the embedding
\(C\hookrightarrow D_J\). Obviously, \(|\sigma(C^\nu)| = i_{D_J}(C,
\partial D_J, Q\cap \partial D_J)\) by definition.

For each point \(p \in \sigma(C^\nu)\) with \(\nu_c(p)\not\in Q\),
\begin{equation}\label{s4:ss3:e6}
\pi_\Gamma^{-1}(p) \in \sigma(\Gamma^\nu)
\end{equation}
by \lemref{lem2}. And for each point \(q \in \sigma(C^\nu)\) with
\(\nu_c(q)\in Q\),
\begin{equation}\label{s4:ss3:e7}
\pi_\Gamma^{-1}(q) \in \sigma(\Gamma^\nu)
\end{equation}
by \lemref{lem3}. Therefore, by Riemann-Hurwitz,
we have
\begin{equation}\label{s4:ss3:e8}
\begin{split}
2 p_a(\Gamma) - 2 + \left|\sigma(\Gamma)\right| &\ge
2 g(\Gamma^\nu) - 2 + \left|\sigma(\Gamma^\nu)\right|\\
&\ge \gamma\left(2g(C^\nu) - 2 + \left|\sigma(C^\nu)\right|\right)\\
&= \gamma \left(2g(C) - 2 + i_{D_J}(C, \partial D_J, Q\cap \partial
D_J)\right)\\
&= \gamma \phi_Q(C).
\end{split}
\end{equation}
This finishes the proof of \eqref{s4:ss3:e2} and hence \thmref{t5}.

\subsection{Proof of \thmref{t8}}\label{s4:ss4}

Let $X\subset\P^3\times\Delta$ be a family of surfaces containing
$D$ where $X_0 = R\cup H_1\cup H_2\cup...\cup H_{d-s}$ with $d-s$
planes $H_1, H_2, ..., H_{d-s}$ in general position and $X$ has only
rational double points at \(Q = X_{sing}\).

Let \(Y\) be a family of curves with the commutative diagram
\eqref{t5:e1}. We want to prove \eqref{t8:e1} holds for $C = Y_t$. By
\thmref{t5}, it suffices to prove
\begin{equation}\label{s4:ss4:e1}
2g_Q^{vir}(\pi(Y_0)) - 2 \ge \epsilon\deg \pi(Y_0)
\end{equation}
with \(\epsilon\) defined as in \eqref{t8:e2}.
Let $\Gamma$ be an irreducible component of $\pi(Y_0)$.

Let $B_i = R\cap H_i$ for $i=1,2,...,d-s$ and $L_{ij} = H_i\cap H_j$
for $1\le i \ne j \le d-s$.

If $\Gamma = B_i$ for some $i$, then
\begin{equation}\label{s4:ss4:e2}
\phi_Q(\Gamma) = 2 g(B_i) - 2 + \sum_{j\ne i} H_j\cdot B_i =
(d-4)\deg\Gamma
\end{equation}
where $\phi_Q(\Gamma)$ is defined in \eqref{defn4:e1}.

If $\Gamma = L_{ij}$ for some $i\ne j$, then
\begin{equation}\label{s4:ss4:e3}
\phi_Q(\Gamma) = 2g(L_{ij}) - 2 + R\cdot L_{ij} + \sum_{k\ne i, j}
H_k \cdot L_{ij} = (d - 4)\deg\Gamma.
\end{equation}

Suppose that $\Gamma\subset R$ and $\Gamma\not\subset H_1\cup H_2\cup
...\cup H_{d-s}$. First, let us assume that $R$ is smooth.

Note that $Q = (D\cup E)\cap (B_1\cup B_2\cup
...\cup B_{d-s})$. Then by \thmref{t2},
\begin{equation}\label{s4:ss4:e4}
\begin{split}
\phi_Q(\Gamma) &= 2g(\Gamma) - 2 + i_R(\Gamma, B_1\cup B_2\cup ...\cup
B_{d-s}, Q)\\
&\ge \min\left(d - 5, d - s -2, \frac{2 g(D) - 2}{\deg D}, \frac{2 g(E)
- 2}{\deg E}\right)\deg\Gamma.
\end{split}
\end{equation}

Suppose that $R$ is not smooth. Obviously, \(s\ge 2\).

Since $E$ is smooth and irreducible,
$R$ is irreducible and smooth in codimension $1$ and hence $B_i$
($1\le i\le d-s$) are irreducible and smooth.
 
Let $\nu: \wt{R}\to R\subset\P^3$
be a minimal desingularization of $R$, $\wt{B}_i =
\nu^{-1} B_i$ for $1\le i\le d-s$, $\wt{Q} = \nu^{-1}(Q)$
and $\wt{D}, \wt{E}, \wt{\Gamma}$ be the proper transforms of $D, E,
\Gamma$, respectively. Instead of applying \thmref{t2} to $(R, \cup
B_i, Q)$, we need to apply it to $(\wt{R}, \cup\wt{B}_i, \wt{Q})$. It
turns out that we only have to show
\begin{equation}\label{s4:ss4:e4.1}
\left(K_{\wt{R}} + (d - s - 1)\wt{B}\right)
\wt{\Gamma} \ge \min(d - 5, d - s - 3)\deg\Gamma
\end{equation}
where $\wt{B}$ is a divisor
representing the divisor class of $\wt{B}_i$, i.e., $\wt{B} = \nu^*(H)$
with $H$ the hyperplane divisor of $\P^3$. Actually, we
claim that \(K_{\wt{R}} + 2 \wt{B}\) is NEF
and it follows that $(K_{\wt{R}} + (d - s - 1)\wt{B})
\wt{\Gamma}\ge (d - s - 3)\deg\Gamma$ and \eqref{s4:ss4:e4.1}. So it
suffices to justify the NEFness of \(K_{\wt{R}} + 2 \wt{B}\).

By Mori's cone theorem, we just have to verify \((K_{\wt{R}} + 2
\wt{B}) G \ge 0\) for all smooth rational curves $G$ such that 
\(-3\le K_{\wt{R}}\cdot G < 0\).

If $K_{\wt{R}}\cdot G = -1$, $G$ is represented by a $-1$ rational
curve. Due to our assumption on $\wt{R}$, $\wt{B} \cdot G > 0$
and hence \((K_{\wt{R}} + 2 \wt{B}) G \ge 0\).

If $K_{\wt{R}}\cdot G = -2$, there is a one-parameter family of curves
in the numerical class of $G$ and hence $G$ is not contracted by $\nu:
\wt{R}\to \P^3$. Therefore, $\wt{B}\cdot G > 0$ and
hence \((K_{\wt{R}} + 2 \wt{B}) G \ge 0\).

If $K_{\wt{R}}\cdot G = -3$, there is a two-parameter family of curves
in the numerical class of $G$ and $\wt{R}$ is rational. So it is
enough to show that \(K_{\wt{R}} + 2 \wt{B}\) is effective. Then the
argument in the proof of \coref{cor0} applies since \((K_{\wt{R}} +
2\wt{B}) \wt{B} = s - 2\ge 0\). Therefore, 
\(K_{\wt{R}} + 2\wt{B}\) is effective and \((K_{\wt{R}} + 2\wt{B})G\ge
0\). Consequently, \(K_{\wt{R}} + 2\wt{B}\) is NEF.

In summary, we have
\begin{equation}\label{s4:ss4:e4.2}
\begin{split}
\phi_Q(\Gamma) &= 2g(\Gamma) - 2 + i_R(\Gamma, B_1\cup B_2\cup ...\cup
B_{d-s}, Q)\\
&\ge \min\left(d - 5, d - s - 3, \frac{2 g(D) - 2}{\deg D}, \frac{2 g(E)
- 2}{\deg E}\right)\deg\Gamma
\end{split}
\end{equation}
for each $\Gamma\subset R$ and $\Gamma\not\subset H_1\cup H_2\cup
...\cup H_{d-s}$.

Suppose that $\Gamma\subset H_i$ for some $i$, say $\Gamma\subset H_1$
and $\Gamma\not\subset R\cup H_2\cup H_3\cup ...\cup H_{d-s}$.
We need to be a little more careful here when we apply \thmref{t2} to
$(H_1, B_1\cup L_{12}\cup L_{13}\cup ...\cup L_{1, d-s}, Q)$. 
Note that $B_1$ is not a general curve of degree $s$ on $H_1\isom \P^2$
since we cannot fix $H_1$ and vary $R$. Hence we should regard $B_1$ as a
fixed curve on $H_1$ while $L_{12}, L_{13}, ..., L_{1, d-s}$
are $d-s-1$ lines in general position on $H_1$.

Note that $Q$ is cut out on $B_1\cup L_{12}\cup L_{13}\cup
...\cup L_{1, d-s}$ by a smooth curve $F\subset H_1$ of
degree $d$. When applying \thmref{t2}, we need to verify the
following:
\begin{equation}\label{s4:ss4:e5}
(K_{H_1} + B_1 + L_{12} + L_{13} + ... + L_{1, d-s-1}) C 
= (d-5) \deg C
\end{equation}
\begin{equation}\label{s4:ss4:e6}
\begin{split}
(L_{12} + L_{13} + ... + L_{1, d-s})C - 2 &= (d-s-1)\deg C - 2\\
& \ge (d-s-3) \deg C
\end{split}
\end{equation}
and
\begin{equation}\label{s4:ss4:e7}
2g(F) - 2 = d(d-3) = (d-3)\deg F.
\end{equation}
Therefore by \thmref{t2}, we have
\begin{equation}\label{s4:ss4:e8}
\begin{split}
\phi_Q(\Gamma) &= 2g(\Gamma) - 2 + i_{H_1}(\Gamma, B_1\cup L_{12}\cup
L_{13}\cup ...\cup L_{1, d-s}, Q)\\
&\ge \min(d-s-3, d-5)\deg\Gamma
\end{split}
\end{equation}
for each \(\Gamma\subset H_1\) and 
\(\Gamma\not\subset R\cup H_2\cup ...\cup H_{d-s}\).

Combining \eqref{s4:ss4:e2}, \eqref{s4:ss4:e3}, \eqref{s4:ss4:e4.2} and
\eqref{s4:ss4:e8}, we obtain \eqref{s4:ss4:e1} and this finishes the
proof of \thmref{t8}.

\subsection{Proof of \thmref{t9}}\label{s4:ss5}

Under the order we choose for \(d_{1j}\) and \(d_{2i}\), we have
\begin{equation}\label{t9:e4}
u_{i+1, j}\le u_{ij}\le u_{i,j+1}
\end{equation}
and by the definition of \(u_{ij}\),
\begin{equation}\label{t9:e5}
u_{ij} + u_{kl} = u_{il} + u_{kj}.
\end{equation}
There is an implicit relation among \(d_{1j}\) and \(d_{2i}\):
\begin{equation}\label{t9:e6}
\sum_{j=1}^{m+1} d_{1j} = \sum_{i=1}^{m} d_{2i}.
\end{equation}
Therefore,
\begin{equation}\label{t9:e7}
d_{1j} = \sum_{k\ne j} u_{kk} + u_{j,m+1}
\end{equation}
for \(j = 1,2,...,m+1\) (let \(u_{m+1,m+1} = 0\)) and
\begin{equation}\label{t9:e8}
d_{2i} = \sum_{k} u_{kk} + u_{i,m+1}
\end{equation}
for \(i = 1,2,...,m\).

More explicitly, \(S\) can be described as the following. It is the
locus given by
\begin{equation}\label{t9:e9}
\rank
\begin{pmatrix}
A_{11} & A_{12} & ... & A_{1m} & A_{1,m+1}\\
A_{21} & A_{22} & ... & A_{2m} & A_{2,m+1}\\
\vdots & \vdots & \ddots & \vdots & \vdots\\
A_{m1} & A_{m2} & ... & A_{mm} & A_{m,m+1}
\end{pmatrix}
< m
\end{equation}
where \(A_{ij}\in H^0(\CO_{\P^4}(u_{ij}))\) is a homogeneous
polynomial of degree \(u_{ij}\) (\(A_{ij} = 0\) if \(u_{ij} < 0\)).
And the homogeneous coordinate ring \(\oplus H^0(I_S(n))\) of \(S\) is
generated by \(F_1, F_2, ..., F_{m+1}\), where \(F_j\) is the minor of
the matrix \((A_{ij})\) obtained by removing its \(j\)-th
column. Obviously, \(\deg F_j = d_{1j}\). Since we assume that
\eqref{t9:e1} is a minimal free resolution and \(S\) is irreducible,
\begin{equation}\label{t9:e10}
u_{ii} > 0 \text{ and } u_{i+1, i}\ge 0.
\end{equation}

Let \(u = u_{m,m+1}\) and choose \(u\) linear forms \(L_1, L_2, ...,
L_u\in H^0(\CO_{\P^4}(1))\).
Let \(X\subset \P^4\times\Delta\) be a family of PCM surfaces given by
\begin{equation}\label{t9:e11}
\rank
\begin{pmatrix}
A_{11} & A_{12} & ... & A_{1m} & L_1 L_2 ... L_u G_1 + t A_{1,m+1}\\
A_{21} & A_{22} & ... & A_{2m} & L_1 L_2 ... L_u G_2 + t A_{2,m+1}\\
\vdots & \vdots & \ddots & \vdots & \vdots\\
A_{m1} & A_{m2} & ... & A_{mm} & L_1 L_2 ... L_u G_m + t A_{m,m+1}
\end{pmatrix}
< m
\end{equation}
where \(G_i\) is a homogeneous polynomial of degree \(u_{i,m+1} - u\)
for \(i=1,2,...,m\). The central fiber of \(X\) consists of \(X_0 =
R\cup H_1\cup H_2 \cup ...\cup H_u\), where \(R\) is a PCM surface of
type \((d_{11} - u, d_{12} - u, ..., d_{1m} - u, d_{1,m+1}, d_{21} -
u, d_{22} - u, ..., d_{2m} - u)\) given by
\begin{equation}\label{t9:e12}
\rank
\begin{pmatrix}
A_{11} & A_{12} & ... & A_{1m} & G_1\\
A_{21} & A_{22} & ... & A_{2m} & G_2\\
\vdots & \vdots & \ddots & \vdots & \vdots\\
A_{m1} & A_{m2} & ... & A_{mm} & G_m
\end{pmatrix}
< m
\end{equation}
and \(H_k\) is a surface given by \(F_{m+1} = L_k = 0\) for
\(k=1,2,...,u\). It is not hard to see that \(X_0\) is of normal
crossing and \(X\) has only rational double points at \(Q =
X_{sing}\). Note that it is possible that \(R = \emptyset\).

Let \(Y\) be a family of curves with the commutative diagram
\eqref{t5:e1}. We want to prove \eqref{t9:e3} holds for \(C = Y_t\). By
\thmref{t5}, it suffices to prove
\begin{equation}\label{t9:e13}
2g_Q^{vir}(\pi(Y_0)) - 2 \ge (u - 5) \deg \pi(Y_0).
\end{equation}
Let $\Gamma$ be an irreducible component of $\pi(Y_0)$. Our goal is to
show
\begin{equation}\label{t9:e14}
\phi_Q(\Gamma) \ge (u-5) \deg \Gamma
\end{equation}
which should imply \eqref{t9:e13} immediately by the definition of
virtual genus.

Let us first check the cases \(\Gamma = R\cap H_i\) or \(\Gamma =
H_i\cap H_j\). If \(\Gamma = R\cap H_i\) for some \(i\), then
\begin{equation}\label{t9:e15}
\begin{split}
\phi_Q(\Gamma) &= 2 g(\Gamma) - 2 + \sum_{k\ne i} H_k\cdot \Gamma\\
&\ge (u-1) \deg \Gamma - 2 \ge (u-3)\deg \Gamma.
\end{split}
\end{equation}
If \(\Gamma = H_i\cap H_j\) for some \(i\ne j\), then
\begin{equation}\label{t9:e16}
\begin{split}
\phi_Q(\Gamma) &\ge 2 g(\Gamma) - 2 + \sum_{k\ne i, j} H_k\cdot \Gamma\\
&\ge (u-2) \deg \Gamma - 2 \ge (u-4)\deg \Gamma.
\end{split}
\end{equation}
In summary, \eqref{t9:e14} holds for \(\Gamma = R\cap H_i\) or \(\Gamma =
H_i\cap H_j\).

Next, let us check the case \(\Gamma \subset R\) and \(\Gamma
\not\subset \cup H_k\). Let \(B_k = R\cap H_k\) and \(B =
\cup B_k\). Basically, we want to apply \thmref{t2} to \((R, B, Q\cap
B)\). There are a few things to check:
\begin{equation}\label{t9:e17}
\left(K_R + B - B_k\right)\Gamma \ge (u - 4) \deg \Gamma.
\end{equation}
for \(1\le k \le u\),
\begin{equation}\label{t9:e18}
B \cdot \Gamma - 2\ge (u-2)\deg \Gamma
\end{equation}
and
\begin{equation}\label{t9:e19}
2 g(C) - 2 \ge (u-3) \deg C
\end{equation}
for all irreducible components \(C\subset F\), where \(F\subset R\) is the
curve such that \(Q\cap B = F\cap B\).

Both \eqref{t9:e17} and \eqref{t9:e18} are more or less obvious; among
the two, \eqref{t9:e17} is due to the fact that \(K_R + 3H\) is NEF
for the hyperplane divisor \(H\). It
remains to verify \eqref{t9:e19}. This requires us to study the curve
\(F\).

By \eqref{t9:e11}, we see that \(F\) is given by
\begin{equation}\label{t9:e20}
\rank
\begin{pmatrix}
A_{11} & A_{12} & ... & A_{1m} & A_{1,m+1} & G_1\\
A_{21} & A_{22} & ... & A_{2m} & A_{2,m+1} & G_2\\
\vdots & \vdots & \ddots & \vdots & \vdots & \vdots\\
A_{m1} & A_{m2} & ... & A_{mm} & A_{m,m+1} & G_m
\end{pmatrix}
< m.
\end{equation}
We claim that such \(F\) is smooth and irreducible for a general choice of
\(A_{ij}\) and \(G_i\) and its Hilbert polynomial (and hence its
degree and arithmetic genus) can be calculated in an inductive way. 

Let us first recall some basic definitions and fix some notations.
For a closed subscheme \(F\subset \P^N\), its Hilbert polynomial \(P_F(l)\)
is the Euler characteristic of the invertible sheaf \(\CO_F(l)\), i.e.,
\begin{equation}\label{t9:e21}
P_F(l) = \chi(\CO_F(l)) = \sum_{k=0}^{\dim F} (-1)^k h^k(\CO_F(l)).
\end{equation}
The leading term of \(P_F(l)\) is \((d/r!) l^r\), where \(r =
\dim F\) and \(d = \deg F\). And the arithmetic genus of \(F\) is
\(p_a(F) = (-1)^r (P_F(0) - 1)\).

Let \(\triangledown\) be the difference operator such that
\(\triangledown f(l) = f(l) - f(l-1)\). Then \(\triangledown
P_F(l)\) is the Hilbert polynomial of a hyperplane section of
\(F\) and \(\triangledown^k P_F(l)\) is the Hilbert polynomial
of an \((N - k)\)-plane section of \(F\).

\begin{prop}\label{prop2}
Let \((u_{ij})_{m\times n}\) (\(m\le n\)) be a matrix of integers satisfying
\eqref{t9:e4}, \eqref{t9:e5} and \eqref{t9:e10}
and \(M = (A_{ij})_{m\times n}\) be a matrix
of homogeneous polynomials, where \(A_{ij}\) is a general member of
\(H^0(\CO_{\P^N}(u_{ij}))\). Let \(F\)
be the locus in \(\P^N\) defined by \(\rank(M) < m\).
Then
\begin{enumerate}
\item \(F\) is irreducible of codimension \(n - m + 1\) in \(\P^N\);
\item \(F_{sing}\) consists of the points such that
\(\rank(M) < m - 1\)
and hence \(F\) is smooth in codimension \(n - m + 2\);
\item its Hilbert polynomial \(P_F(l)\) can be computed inductively in
the following way.
\end{enumerate}

Assume that \(m < n\).
Pick \(m \le \alpha \le n\) and let \(M_\alpha\) be the matrix
obtained from \(M\) by removing the \(\alpha\)-th column of \(M\) and
\(M_{m\alpha}\) be the matrix obtained from \(M\) by removing the
\(m\)-th row and \(\alpha\)-th column. And let \(F_\alpha\) and
\(F_{m\alpha}\) be the corresponding subvarieties of \(\P^N\) given by
\(\rank(M_\alpha) < m\) and \(\rank(M_{m\alpha}) < m - 1\). Then
\begin{equation}\label{prop2:e1}
P_F(l) = \sum_{k\ge 0} (-1)^k \left(\binom{u}{k+1}
\triangledown^{k+1} P_{F_\alpha} (l) +
\binom{u}{k} \triangledown^{k} P_{F_{m\alpha}} (l)\right),
\end{equation}
where \(u = u_{m\alpha}\).
\end{prop}

\begin{proof}
The first two statements follow from a standard ``projection''
argument. Let 
\begin{equation}\label{prop2:e2}
W\subset \P\left(\prod_{i,j} H^0(\CO_{\P^N}(u_{ij}))\right)\times \P^N
\end{equation} 
be the incidence correspondence given by \(\rank(A_{ij}) < m\),
where we regard \(A_{ij}\) as a homogeneous polynomial with generic
coefficients. Obviously, \(F\) is a general fiber of \(W\) when it is
projected onto the first factor \(\P\left(\prod_{i,j}
H^0(\CO_{\P^N}(u_{ij}))\right)\). To see that \(F\) has the stated
properties, it suffices to justify that \(W\) has the same properties,
i.e., \(W\) is irreducible of codimension \(n-m+1\) in the total space
and \(W_{sing}\) consists of the points such that \(\rank(A_{ij}) <
m-1\). This can be shown by projecting \(W\) to the second factor
\(\P^N\). The general fiber of \(W\) over a point in \(\P^N\) is
(almost) a generic determinantal variety, which is classically known to have
the stated properties. We will let the readers to fill out the details.

Choose linear forms \(L_1, L_2, ..., L_u\) and specialize
the \(\alpha\)-th column of \(M\) to \(A_{i\alpha} = L_1 L_2 ... L_u
G_i\). Since Hilbert polynomials are invariant under deformations,
such degeneration does not change \(P_F(l)\). Let \(H_k\) be the
hyperplane given by \(L_k = 0\). Then 
\(F = \cup_{k=1}^u (H_k\cap F_\alpha) \cup Z\), where \(Z\) is the
subvariety of \(\P^N\) given by
\begin{equation}\label{prop2:e3}
\rank
\begin{pmatrix}
A_{11} & ... & A_{1m} & ... & G_1 & ... & A_{1n}\\
A_{21} & ... & A_{2m} & ... & G_2 & ... & A_{2n}\\
\vdots & \ddots & \vdots & \ddots & \vdots & \ddots & \vdots\\
A_{m1} & ... & A_{mm} & ... & G_m & ... & A_{mn}
\end{pmatrix}
< m
\end{equation}
where \(G_i\)'s occupy the \(\alpha\)-th column of the matrix and the
rest of the matrix are the same as \(M\). Observe
that \(\deg G_m = 0\) and hence we may apply row operations to the
matrix to eliminate its \(m\)-th row and \(\alpha\)-th
column. Therefore, it is easy to see that \(Z\) and \(F_{m\alpha}\)
are deformationally equivalent inside \(\P^N\). Especially they have
the same numerical properties.

For a general choices of \(L_k\), \(F\) is of normal crossing with
components \(H_k\cap F_\alpha\) and \(Z\). Apply
the inclusion-exclusion principle to \(P_F(l)\) and the recursion
formula \eqref{prop2:e1} follows.
\end{proof}

By the above proposition, \(F\) is a smooth
irreducible curve. Therefore, in order to prove \eqref{t9:e19}, it suffices
to show that
\begin{equation}\label{t9:e22}
-2 \chi(\CO_F) = 2p_a(F) - 2 \ge (u-3) \deg F
\end{equation}
where \(p_a(F)\) and \(\deg F\) can be computed, in principle, by
\eqref{prop2:e1}. However, it takes a substantial amount of
computation to write down
\(p_a(F)\) and \(\deg F\) explicitly in terms of \(u_{ij}\). Nor is this
necessary for our cause. All we need to do is to prove \eqref{t9:e22}
by an induction via \eqref{prop2:e1}.

\begin{prop}\label{prop3}
Let \((u_{ij})_{m\times n}\), \(M = (A_{ij})_{m\times n}\) and \(F\)
be defined as in \propref{prop2}. Here we let \(N = n - m + 2\) and
hence \(F\) is a curve. Then \eqref{t9:e22} holds for \(u = u_{mn}\).
\end{prop}

\begin{proof}
We do induction on \(\min(n - m, m)\). If \(m = n\) or \(m = 1\),
\(F\) is a complete intersection in \(\P^N\) and \eqref{t9:e22}
follows easily from adjunction.

Suppose that \(n > m > 1\). Let \(\alpha = m\) in \eqref{prop2:e1} and
we have
\begin{equation}\label{prop3:e1}
-\chi(\CO_F) \ge -u_{mm} \chi(\CO_{F_{m}\cap H}) - \chi(\CO_{F_{mm}})
\end{equation}
and
\begin{equation}\label{prop3:e2}
\deg F = u_{mm} \deg (F_{m}\cap H) + \deg F_{mm}
\end{equation}
where \(F_m\) and \(F_{mm}\) are defined as in \propref{prop2} and
\(H\) is a hyperplane in \(\P^N\). 

By induction hypothesis, \(- 2\chi(\CO_{F_{m}\cap H}) \ge (u-3) \deg
(F_{m}\cap H)\) and \(- 2\chi(\CO_{F_{mm}}) \ge (u-3) \deg
F_{mm}\). Therefore, \eqref{t9:e22} follows from \eqref{prop3:e1} and
\eqref{prop2:e2}.
\end{proof}

Combining \eqref{t9:e17}-\eqref{t9:e19} and applying \thmref{t2} to
\((R, B, Q\cap B)\), we see that \eqref{t9:e14} holds for all
\(\Gamma\subset R\) and \(\Gamma\not\subset \cup H_k\).

Finally, we need to check the case \(\Gamma\subset H_i\)
and \(\Gamma \not\subset \cup_{k\ne i} H_k\cup R\) for some \(i\).
Without the loss of generality, let us assume that \(i = u\),
\(\Gamma \subset H = H_u\), \(B_0 = R\cap H\), \(B_k = H_k\cap H\) for
\(1\le k < u\) and \(B = B_1\cup B_2 \cup ... \cup B_{u-1}\). Again we
want to apply \thmref{t2} to \((H, B, Q\cap B)\). Here we regard
\(B_0\) as a fixed curve while \(B_k\) varies in the BPF linear system
cut out by the hyperplane divisor for \(1\le k < u\).
So we need to verify the following:
\begin{equation}\label{t9:e23}
(K_H + B - B_k)\Gamma \ge (u - 5) \deg \Gamma
\end{equation}
for \(1\le k < u\)
\begin{equation}\label{t9:e24}
(B - B_0)\Gamma - 2 \ge (u - 3) \deg \Gamma
\end{equation}
and
\begin{equation}\label{t9:e25}
2 g(C) - 2 \ge (u-3) \deg C
\end{equation}
for all irreducible components \(C\subset F\), where \(F\subset H\) is
the curve such that \(Q\cap B = F\cap B\).

Again, both \eqref{t9:e23} and \eqref{t9:e24} are obvious; among the
two, \eqref{t9:e23} follows from Mori's cone theorem. To see
\eqref{t9:e25}, we observe that \(F\) is a curve in \(\{ L_u = 0\}\isom
\P^3\) defined by \eqref{t9:e9} with \(A_{ij}\) restricted to the
hyperplane \(\{ L_u
= 0\}\). So by \propref{prop2} and \ref{prop3}, \(F\) is irreducible
and smooth and \eqref{t9:e22} holds. Therefore, \eqref{t9:e25}
follows and this finishes the proof of \eqref{t9:e14} and hence the
proof of \thmref{t9}.

\end{document}